# SECOND-ORDER FLUCTUATIONS AND CURRENT ACROSS CHARACTERISTIC FOR A ONE-DIMENSIONAL GROWTH MODEL OF INDEPENDENT RANDOM WALKS[1]

By Timo Seppäläinen

*University of Wisconsin*

Fluctuations from a hydrodynamic limit of a one-dimensional asymmetric system come at two levels. On the central limit scale $n^{1/2}$ one sees initial fluctuations transported along characteristics and no dynamical noise. The second order of fluctuations comes from the particle current across the characteristic. For a system made up of independent random walks we show that the second-order fluctuations appear at scale $n^{1/4}$ and converge to a certain self-similar Gaussian process. If the system is in equilibrium, this limiting process specializes to fractional Brownian motion with Hurst parameter $1/4$. This contrasts with asymmetric exclusion and Hammersley's process whose second-order fluctuations appear at scale $n^{1/3}$, as has been discovered through related combinatorial growth models.

**1. Introduction.** An interface model defined in terms of a height function on an integer lattice is a stochastic process $\sigma_t = \{\sigma_t(x) : x \in \mathbf{Z}^d\}$, where the value $\sigma_t(x)$ is interpreted as the height of the interface over site $x$. The random variables $\sigma_t(x)$ move up and down according to random rates whose momentary values depend on the height values in some neighborhood around site $x$. A hydrodynamic scaling limit is a type of law of large numbers for these systems. The conventional statement is that under a suitable scaling of space and time, the entire space-time random evolution $\{\sigma_t(x) : x \in \mathbf{Z}^d, t \geq 0\}$ converges to the solution of a differential equation.

When the system is asymmetric, in the sense that there is an average drift either up or down, the typical law of large numbers for the system is of the following form. The result is for a sequence of processes $\sigma_t^n$, where $n = 1, 2, 3, \ldots$ is the scaling parameter. The statement is that as $n \to \infty$, the

Received August 2003; revised April 2004.
[1]Supported in part by NSF Grant DMS-01-26775.

random position $n^{-1}\sigma^n_{nt}([nx])$ converges to a function $u(x,t)$ that solves a Hamilton–Jacobi equation

$$u_t + f(\nabla u) = 0. \tag{1.1}$$

In the macroscopic description, $f$ gives the local velocity of the height as a function of the local gradient. A necessary hypothesis for this type of law of large numbers is that as $n \to \infty$, the sequence of scaled initial states $n^{-1}\sigma^n_0([ny])$ in some sense converges to a function $u_0(y)$. The function $u_0$ then serves as the initial data for the equation $u_t + f(\nabla u) = 0$. The equation and $u_0$ uniquely determine $u(x,t)$ at all later times $t$. Depending on the situation, there might be additional assumptions on the distributions of $\sigma^n_0$. Examples of results for various models can be found in [12, 14, 16, 17, 19]. For general accounts of hydrodynamic limits, we refer to the monographs [11] and [20], and to the lectures [4] and [21].

The fluctuation question for asymmetric systems has so far found answers only in the one-dimensional situation, where the following picture has emerged. Suppose the initial conditions satisfy a central limit theorem of the type

$$\frac{\sigma^n_0([ny]) - nu_0(y)}{\sqrt{n}} \to \zeta_0(y), \qquad y \in \mathbf{R},$$

with a continuous limiting process $\zeta_0$. This situation arises naturally when the initial increments $\sigma^n_0(x) - \sigma^n_0(x-1)$ are independent with slowly varying bounded means and variances. Then at later times a weak limit

$$\frac{\sigma^n_{nt}([nx]) - nu(x,t)}{\sqrt{n}} \to \zeta(x,t) \tag{1.2}$$

holds. The process $\zeta(x,t)$ is a deterministic function of the initial process $\zeta_0$. More specifically, the value $\zeta(x,t)$ is determined by the values $\zeta_0(y)$, such that a generalized characteristic of the p.d.e. (1.1) emanating at $(y,0)$ reaches $(x,t)$. Qualitatively, a crucial feature is that there is no dynamical noise visible at this scale $n^{1/2}$. These types of results have been proved for the exclusion process under various hypotheses [7, 8, 13] and for Hammersley's process [18].

The motivation of the present paper is to describe fluctuations that lie "beyond" the trivial fluctuations transported by the characteristics that appear in (1.2). This second level of fluctuations appears when the first-order fluctuations of (1.2) are suitably subtracted off or when the initial conditions are deterministic. In the asymmetric exclusion and Hammersley settings, one should find some kind of fluctuations at the $n^{1/3}$ scale. The results of Baik–Deift–Johansson [2] and Johansson [10] can be interpreted as fluctuation results for Hammersley's process and the exclusion process with special deterministic initial configurations. So for exclusion and Hammersley's process, the aim would be to generalize those $n^{1/3}$ results to other initial conditions. What is special about exclusion and Hammersley's process is that these "second-order" fluctuations arise through natural growth models that are amenable to the powerful combinatorial and analytic machinery of [2] and its descendants.

In the present paper we start another direction, the investigation of these phenomena in other models besides exclusion and Hammersley's process. Two questions arise. Do the trivial $n^{1/2}$ fluctuations transported by characteristics appear in other asymmetric models? What then would be the second-order fluctuations, especially if there is no combinatorial growth model present that would lead to the $n^{1/3}$ fluctuations and the random matrix connections?

A simple model is one where the increments of the height function come from independent random walks. This case we analyze in the present paper. The initial increments are taken independently with slowly varying means and variances. Exactly as for exclusion and Hammersley, on the central limit scale $n^{1/2}$ we have the initial fluctuations transported by characteristics. Then we find the next order of fluctuations on the scale $n^{1/4}$. In the limit these fluctuations are described by a certain Gaussian process $Z(\bar{y}, t)$, where $(\bar{y}, 0)$ is the initial point of the (now unique) characteristic that reaches $(x, t)$. The covariance of $Z(\bar{y}, \cdot)$ is determined by the mean and variance of the initial increments around the macroscopic point $\bar{y}$. Imprecisely speaking, the random height expands as

$$\sigma_{nt}^n([nx]) \approx nu(x,t) + n^{1/2}\zeta(x,t) + n^{1/4}Z(\bar{y},t).$$

The processes $Z(\bar{y}, \cdot)$ are independent for distinct initial points $\bar{y}$.

In the special case when the height increments are in equilibrium, for a fixed $\bar{y}$, the process $Z(t) = Z(\bar{y}, t)$ specializes to fractional Brownian motion with Hurst parameter $\frac{1}{4}$. These second-order fluctuations turn out to be the same as the particle current across a characteristic of the macroscopic equation. Hence, the juxtaposition in the title of the paper.

In the next section we describe the random walk model, its hydrodynamic limit and the two levels of fluctuations. For the sake of comparison, we include a brief section on the fluctuation picture of Hammersley's process. We show that for Hammersley's process the second-order fluctuations are, at most, of order $n^{1/3} \log n$. This bound is valid also for shock locations. The last two sections prove the theorems. In the proofs, $C$, $C_0$, $C_1, \ldots$ denote constants whose actual values may change from line to line. The set of natural numbers is $\mathbf{N} = \{1, 2, 3, \ldots\}$.

**2. The random walk model.** We consider a model of an interface whose height differences between neighboring sites are defined by independent random walks on $\mathbf{Z}$. The state of the system at time $t$ is a height function

$\sigma_t : \mathbf{Z} \to \mathbf{Z}$. The height function is nondecreasing in space, so the increments

$$\eta_t(x) = \sigma_t(x) - \sigma_t(x-1) \tag{2.1}$$

are nonnegative integers.

The randomly evolving height function is constructed as follows. Let $\{X_i(t) : i \in \mathcal{I}\}$ be a countable collection of independent continuous-time random walks on $\mathbf{Z}$. The jump rates of the random walks are given by a probability kernel $\{p(x) : x \in \mathbf{Z}\}$. In other words, the assumption on $p(x)$ is

$$\sum_{x \in \mathbf{Z}} p(x) = 1$$

and the common transition probability of the random walks is

$$P[X_i(s+t) = y | X_i(s) = x] = p_t(x, y) = \sum_{k=0}^{\infty} \frac{e^{-t} t^k}{k!} p^{(k)}(y-x),$$

where

$$p^{(k)}(z) = \sum_{x_1 + x_2 + \cdots + x_k = z} p(x_1) p(x_2) \cdots p(x_k)$$

is the $k$-fold convolution of the kernel $p(x)$.

Given an initial height function $\sigma_0 = \{\sigma_0(x) : x \in \mathbf{Z}\}$ defined on some probability space, define the initial increments $\eta_0(x) = \sigma_0(x) - \sigma_0(x-1)$. Choose the initial positions of the random walks so that site $x$ contains $\eta_0(x)$ particles:

$$\sum_{i \in \mathcal{I}} \mathbf{1}\{X_i(0) = x\} = \eta_0(x). \tag{2.2}$$

Once the initial points $\{X_i(0) : i \in \mathcal{I}\}$ have been specified, define the subsequent evolutions $\{X_i(t) - X_i(0) : i \in \mathcal{I}, t \geq 0\}$ as an i.i.d. collection of random walks on this same probability space, independent of $\sigma_0$. Define the current $J_t(x)$ as the (net) number of particles that have moved from $(-\infty, x]$ to $[x+1, \infty)$ in time interval $[0, t]$,

$$J_t(x) = \sum_i \mathbf{1}\{X_i(0) \leq x < X_i(t)\} - \sum_i \mathbf{1}\{X_i(t) \leq x < X_i(0)\}.$$

The height function at time $t$ is then defined by

$$\sigma_t(x) = \sigma_0(x) - J_t(x). \tag{2.3}$$

In other words, the interface height at $x$ moves *down* one step with every particle that jumps from $(-\infty, x]$ to $[x+1, \infty)$, and *up* one step with every particle that does the opposite. The increment variables $\eta_t(x)$ defined by

(2.1) also serve as the occupation variables of the random walks: from (2.1)–(2.3) one can derive

$$\eta_t(x) = \sum_i \mathbf{1}\{X_i(t) = x\}. \tag{2.4}$$

We can describe the evolution of $\sigma_t$ directly in terms of the rates, without reference to the random walks. Given a height function $\sigma$ (a nondecreasing function $\mathbf{Z} \to \mathbf{Z}$), and $x, \ell \in \mathbf{Z}$, define a new height function $\sigma^{x,\ell}$ by

$$\sigma^{x,\ell}(y) = \begin{cases} \begin{cases} \sigma(y) - 1, & x \leq y \leq x + \ell - 1, \\ \sigma(y), & \text{if } y < x \text{ or } y \geq x + \ell, \end{cases} & \text{if } \ell \geq 0; \\ \begin{cases} \sigma(y) + 1, & x + \ell \leq y \leq x - 1, \\ \sigma(y), & \text{if } y < x + \ell \text{ or } y \geq x, \end{cases} & \text{if } \ell < 0. \end{cases} \tag{2.5}$$

The dynamical rule for the height process is this: if the current state is $\sigma$, then for each $x \in \mathbf{Z}$ and $\ell \in \mathbf{Z}$, at rate $p(\ell)(\sigma(x) - \sigma(x-1))$, the process jumps from $\sigma$ to $\sigma^{x,\ell}$. If $\ell = 0$, there is actually no change: $\sigma^{x,0} = \sigma$.

Now some assumptions. We give them in three groups, first the assumption on the kernel $p(x)$ and then the assumptions on the sequence of initial height functions $\{\sigma_0^n\}$ that determine the hydrodynamic limit setting.

ASSUMPTION A. For the random walk kernel, we assume that, for some $\delta > 0$,

$$\sum_{x \in \mathbf{Z}} e^{\theta x} p(x) < \infty \qquad \text{for } |\theta| \leq \delta. \tag{2.6}$$

The purpose of Assumption A is to enable us to use standard large deviation bounds on the random walks.

Assume given a sequence of initial height functions $\sigma_0^n$, random or deterministic, defined on some probability space. The positive integer parameter $n$ will tend to $\infty$ in the results. Define initial occupation variables $\eta_0^n(x) = \sigma_0^n(x) - \sigma_0^n(x-1)$.

ASSUMPTION B. Assume that for some nondecreasing $C^1$ function $u_0$ on $\mathbf{R}$, and all $y \in \mathbf{R}$,

$$\lim_{n \to \infty} n^{-1} \sigma_0^n([ny]) = u_0(y) \qquad \text{in probability.} \tag{2.7}$$

Assumption B is for the hydrodynamic limit. For the fluctuation results we need stricter control of the initial conditions, as in the next assumption.

ASSUMPTION C. For each $n$, the initial occupation variables $\{\eta_0^n(x) : x \in \mathbf{Z}\}$ are independent, with a uniformly bounded sixth moment:

$$\sup_{n \in \mathbf{N}, x \in \mathbf{Z}} E[\eta_0^n(x)^6] < \infty. \tag{2.8}$$

Let
$$\rho_0^n(x) = E\eta_0^n(x) \quad \text{and} \quad v_0^n(x) = \text{Var}[\eta_0^n(x)]$$
be the mean and variance of the initial occupation variable $\eta_0^n(x)$, $x \in \mathbf{Z}$. Let $u_0$ be the function specified in Assumption B, set $\rho_0 = u_0'$, and let $v_0$ be another given nonnegative function on $\mathbf{R}$. Assume both $\rho_0$ and $v_0$ are bounded. The means $\rho_0^n(x)$ and variances $v_0^n(x)$ approximate the functions $\rho_0$ and $v_0$ in the following precise sense:

For each $y \in \mathbf{R}$, there exist positive integers $L = L(n)$ such that $n^{-1/4}L(n) \to 0$, and for any finite constant $A$,

$$(2.9) \quad \lim_{n\to\infty} \sup_{|m| \leq A\sqrt{n\log n}} n^{1/4} \left| \frac{1}{L(n)} \sum_{j=1}^{L(n)} \rho_0^n([ny] + m + j) - \rho_0(y) \right| = 0.$$

The same assumption holds when $\rho_0^n$ and $\rho_0$ are replaced by $v_0^n$ and $v_0$.

Throughout the paper, $L = L(n)$ denotes the quantity specified in the assumption above. The awkwardly complicated assumption (2.9) is made to accommodate both random and deterministic initial conditions. Given a function $\rho_0$, the expectation of a random $\eta_0^n(x)$ can, of course, agree exactly with $\rho_0(\frac{x}{n})$, but a deterministic $\eta_0^n(x)$ cannot unless $\rho_0(\frac{x}{n})$ is integer-valued. Here are two basic examples of initial conditions that satisfy Assumptions B and C.

EXAMPLE 2.1 (*Random initial conditions*). The functions $\rho_0 = u_0'$ and $v_0$ are bounded, nonnegative and satisfy a local Hölder property: for each bounded interval $[a,b]$, there exist $\beta = \beta(a,b) > 1/2$ and $C = C(a,b) < \infty$ such that

$$(2.10) \quad |\rho_0(x) - \rho_0(y)| + |v_0(x) - v_0(y)| \leq C|x-y|^\beta \quad \text{for all } x, y \in [a,b].$$

For each $n$, let $\{\eta_0^n(x) : x \in \mathbf{Z}\}$ be independent, satisfy assumption (2.8) and have

$$(2.11) \quad E\eta_0^n(x) = \rho_0\left(\frac{x}{n}\right) \quad \text{and} \quad \text{Var}[\eta_0^n(x)] = v_0\left(\frac{x}{n}\right).$$

Additionally, the variables $\{\sigma_0^n(0)\}$ are chosen so that (2.7) holds for $y = 0$. Then (2.7) is satisfied for all $y \in \mathbf{R}$.

EXAMPLE 2.2 (*Deterministic initial conditions*). $\rho_0 = u_0'$ is bounded, nonnegative and satisfies the Hölder condition (2.10) with $\beta = \beta(a,b) > 1/2$, and $v_0$ is identically zero. Define deterministic initial occupation variables by

$$(2.12) \quad \eta_0^n(m) = \left[nu_0\left(\frac{m}{n}\right)\right] - \left[nu_0\left(\frac{m-1}{n}\right)\right],$$

where, as throughout the paper, $[x] = \max\{k \in \mathbf{Z} : k \leq x\}$ denotes the integer part of a real $x$. The function $u_0$ is nondecreasing so each $\eta_0^n(m)$ is a nonnegative integer. And finally, the variables $\{\sigma_0^n(0)\}$ can be random or deterministic, but must satisfy (2.7) for $y = 0$.

We set the stage with the hydrodynamic limit and the first-order fluctuations. The first two moments of the random walk kernel appear in various parts of the results. We denote these by

$$b = \sum_x x p(x) \quad \text{and} \quad \kappa_2 = \sum_x x^2 p(x).$$

Let

(2.13) $$u(x,t) = u_0(x - bt).$$

It is the solution of the linear transport equation

(2.14) $$u_t + b u_x = 0, \qquad u(x,0) = u_0(x).$$

THEOREM 2.1. *Assume Assumptions A–C. Then, for each $(x,t) \in \mathbf{R} \times [0,\infty)$,*

(2.15) $$\lim_{n \to \infty} n^{-1} \sigma_{nt}^n([nx]) = u(x,t) \qquad \text{in probability.}$$

*Given $(x,t)$, let $\bar{y} = x - bt$. Then*

(2.16) $$\lim_{n \to \infty} \left\{ \frac{\sigma_{nt}^n([nx]) - nu(x,t)}{\sqrt{n}} - \frac{\sigma_0^n([n\bar{y}]) - nu_0(\bar{y})}{\sqrt{n}} \right\} = 0 \qquad \text{in probability.}$$

*Furthermore, assume the specific situation of Example* 2.1 *and the normalization $u_0(0) = \sigma_0^n(0) = 0$ for all $n$. Then it is possible to construct the processes $\sigma_t^n$ on the same probability space with a two-sided standard Brownian motion $B(\cdot)$ such that these limits hold in probability, for all $(x,t)$, $\bar{y} = x - bt$:*

(2.17) $$\lim_{n \to \infty} \frac{\sigma_{nt}^n([nx]) - nu(x,t)}{\sqrt{n}} = B\left( \int_0^{\bar{y}} v_0(s) \, ds \right).$$

A two-sided Brownian motion $B(\cdot)$ is constructed by taking two independent standard Brownian motions $B_1$ and $B_2$ on $[0,\infty)$ and setting

$$B(t) = \begin{cases} B_1(t), & t \geq 0, \\ -B_2(-t), & t < 0. \end{cases}$$

The integral $\int_0^{\bar{y}} v_0(s) \, ds$ in (2.17) has to be interpreted with a sign, in other words,

$$\int_0^{\bar{y}} v_0(s) \, ds = -\int_{\bar{y}}^0 v_0(s) \, ds \qquad \text{for } \bar{y} < 0.$$

The characteristics of (2.14) are straight lines with slope $b$, so $x = \bar{y} + bt$ is the characteristic starting at $(\bar{y}, 0)$. Limits (2.16) and (2.17) say that on the central limit scale $n^{1/2}$, fluctuations from the hydrodynamic limit consist of initial fluctuations rigidly transported along the characteristics, without any contribution from dynamical noise.

The purpose of this paper is to describe the "second-order" fluctuations that appear beyond the trivial fluctuations of Theorem 2.1. Fix $\bar{y} \in \mathbf{R}$. Let

$$(2.18) \qquad Y_n(t) = \sigma_{nt}^n([n\bar{y}] + [nbt]) - \sigma_0^n([n\bar{y}]).$$

Since $u(\bar{y} + bt, t) = u_0(\bar{y})$, also

$$(2.19) \quad \sigma_{nt}^n([n\bar{y}] + [nbt]) - nu(\bar{y} + bt, t) = \sigma_0^n([n\bar{y}]) - nu_0(\bar{y}) + Y_n(t).$$

So $Y_n(t)$ represents the difference between the fluctuation experienced by the process at space-time point $(\bar{y} + bt, t)$ and the fluctuation at the initial point $(\bar{y}, 0)$ of the characteristic. $n^{-1/2} Y_n(t)$ is exactly the difference that appeared in (2.16).

There is another way to look at $Y_n(t)$, directly in terms of the particles. Write $J_t^n$ for the current of process $\sigma_t^n$ and $X_i^n(t)$ for the random walks in the construction of $\sigma_t^n$:

$$\begin{aligned}
Y_n(t) &= \sigma_{nt}^n([n\bar{y}] + [nbt]) - \sigma_0^n([n\bar{y}]) \\
&= \sigma_0^n([n\bar{y}] + [nbt]) - J_{nt}^n([n\bar{y}] + [nbt]) - \sigma_0^n([n\bar{y}]) \\
&= \sum_{m=[n\bar{y}]+1}^{[n\bar{y}]+[nbt]} \eta_0^n(m) - J_{nt}^n([n\bar{y}] + [nbt]).
\end{aligned}$$

Switching back to the random walks and cancelling gives

$$(2.20) \qquad \begin{aligned} Y_n(t) &= \sum_i \mathbf{1}\{X_i^n(0) \geq [n\bar{y}] + 1, X_i^n(nt) \leq [n\bar{y}] + [nbt]\} \\ &\quad - \sum_i \mathbf{1}\{X_i^n(0) \leq [n\bar{y}], X_i^n(nt) > [n\bar{y}] + [nbt]\}. \end{aligned}$$

This counts the net number of particles that have moved from the right side of the characteristic to the left side during time interval $[0, nt]$. In other words, $Y_n(t)$ also represents the negative of the current across the characteristic.

Our main theorem is the distributional limit of $Y_n$. Assumption B is not relevant for the limit of $Y_n$. The previous paragraph already showed that even though $Y_n(t)$ was defined in terms of the height functions in (2.18), it is actually determined by the increment process $\eta_t^n$. We included Assumption B in the earlier discussion only to give the complete hydrodynamic picture. Also, the approximation assumption (2.9) is needed only for the particular $\bar{y}$ that appears in the definition of $Y_n$.

THEOREM 2.2. *Fix $\bar{y} \in \mathbf{R}$ and define $Y_n(t)$ as in* (2.18). *Assume Assumptions* A *and* C. *Then as $n \to \infty$, the process $n^{-1/4} Y_n(\cdot)$ converges weakly on the space $D_\mathbf{R}[0,\infty)$ to the mean-zero Gaussian process $Z(\cdot)$ with covariance*

$$(2.21) \qquad EZ(s)Z(t) = \sqrt{\frac{\kappa_2}{2\pi}} \{ \rho_0(\bar{y})(\sqrt{s+t} - \sqrt{s \vee t - s \wedge t}) \\ + v_0(\bar{y})(\sqrt{s} + \sqrt{t} - \sqrt{s+t}) \}.$$

For the increment process $\eta_t = \{\eta_t(m) : m \in \mathbf{Z}\}$, i.i.d. Poisson distributions are equilibrium distributions. If the mean of the Poisson is $\rho$, then $\rho_0(x) = v_0(x) = \rho$ for all $x \in \mathbf{R}$. The covariance in (2.21) then simplifies to

$$(2.22) \qquad EZ(s)Z(t) = \rho \sqrt{\frac{\kappa_2}{2\pi}} (\sqrt{s} + \sqrt{t} - \sqrt{t-s}) \qquad \text{for } s \leq t.$$

This is the covariance of fractional Brownian motion with Hurst parameter $H = \frac{1}{4}$, normalized by $EZ(1)^2 = \rho\sqrt{2\kappa_2/\pi}$. *Standard* fractional Brownian motion would have $EZ(1)^2 = 1$. Let us state this special case as a corollary.

COROLLARY 2.1. *Let $\eta_t$ be an equilibrium process whose occupation variables $\{\eta_t(x) : x \in \mathbf{Z}\}$ are i.i.d. Poisson with mean $\rho$. Let $\eta_t^n = \eta_t$ for each $n$. Then $n^{-1/4} Y_n(\cdot)$ converges weakly on the space $D_\mathbf{R}[0,\infty)$ to fractional Brownian motion with covariance given in* (2.22).

Another special case worth highlighting is that of deterministic initial height functions described in Example 2.2. In that case the limit in (2.17) is zero because the second fraction in (2.16) vanishes in the limit. (2.18) gives the actual fluctuations from the hydrodynamic limit because $\sigma_0^n([n\bar{y}])$ is deterministic. We also omit the short derivation of this corollary from Theorem 2.2.

COROLLARY 2.2. *Let $u_0$ be a nondecreasing $C^1$ function on $\mathbf{R}$ with bounded derivative $\rho_0$. Assume $\rho_0$ satisfies the local Hölder condition* (2.10). *Define deterministic initial height functions by $\sigma_0^n(m) = [nu_0(\frac{m}{n})]$, $m \in \mathbf{Z}$. The function $v_0$ is now zero. Fix $\bar{y} \in \mathbf{R}$. Then the process*

$$\left\{ \frac{\sigma_{nt}^n([n\bar{y}] + [nbt]) - nu(\bar{y} + bt, t)}{n^{1/4}} : t \geq 0 \right\}$$

*converges weakly on $D_\mathbf{R}[0,\infty)$ to the mean zero Gaussian process $Z(\cdot)$ with covariance*

$$EZ(s)Z(t) = \rho_0(\bar{y})\sqrt{\frac{\kappa_2}{2\pi}}(\sqrt{s+t} - \sqrt{t-s}) \qquad \text{for } s \leq t.$$

Fractional Brownian motion has stationary increments, but the general process $Z(t)$ with covariance (2.21) does not unless $\rho_0(\bar{y}) = v_0(\bar{y})$. One can check that for a fixed $h > 0$, $E[(Z(t+h) - Z(t))^2]$ is strictly decreasing with $t$ if $v_0(\bar{y}) > \rho_0(\bar{y})$, and strictly increasing if $v_0(\bar{y}) < \rho_0(\bar{y})$. A bound

$$E[(Z(t) - Z(s))^2] \leq C(t-s)^{1/2}$$

is valid for all $0 \leq s < t$. Since the increment is mean-zero Gaussian, it follows from Kolmogorov's criterion that the process $Z$ has continuous paths.

The process $Z$ is self-similar with index $\frac{1}{4}$, which means that $\{Z(at) : t \geq 0\} \stackrel{d}{=} \{a^{1/4} Z(t) : t \geq 0\}$, as is immediate from the form of the covariance.

Next we address the joint distribution of processes $Y_n(\cdot)$ from several initial points $\bar{y}$. Write $Y_n(\bar{y}, t)$ for the random variable defined by (2.18) or, equivalently, by (2.20) to display its dependence on $\bar{y}$. Write $Z(\bar{y}, t)$ for the Gaussian process with covariance given in (2.21).

THEOREM 2.3. *Assume Assumptions* A *and* C. *Let* $\bar{y}_1 < \bar{y}_2 < \cdots < \bar{y}_k$ *be points on* $\mathbf{R}$. *Then as* $n \to \infty$, *the joint process* $n^{-1/4}(Y_n(\bar{y}_1, \cdot), Y_n(\bar{y}_2, \cdot), \ldots, Y_n(\bar{y}_k, \cdot))$ *converges in distribution on the space* $D_{\mathbf{R}^k}[0, \infty)$ *to a vector* $(Z(\bar{y}_1, \cdot), Z(\bar{y}_2, \cdot), \ldots, Z(\bar{y}_k, \cdot))$ *of independent* $D_{\mathbf{R}}[0, \infty)$-*valued components.*

*Remark about mean-zero random walks.* We have made no assumption on the mean $b$ of the random walk. The results are true also for $b = 0$. However, in this case the convergence of $Y_n$ does not relate to the hydrodynamic limit in the same way because Theorem 2.1 is not the correct limit. The relevant hydrodynamic limit takes place on the time scale $n^2 t$ and the limiting evolution is governed by the heat equation. For $b = 0$, Theorem 2.1 is completely trivial because $u(x, t) = u_0(x)$ and $\bar{y} = x$.

*Remark about fractional Brownian motion with* $H = \frac{1}{4}$. There is, of course, a result for Brownian motion that corresponds to the random walk result of Theorem 2.2. We state here the equilibrium version. Let $\lambda > 0$. Let $\{B_i(t) : i \in \mathcal{I}\}$ be a countable collection of independent standard Brownian motions on $\mathbf{R}$ whose initial locations (and, consequently, the locations at any fixed time) are those of a homogeneous, rate $\lambda$ Poisson point process on $\mathbf{R}$. Fix $y \in \mathbf{R}$. Let

$$Y_\lambda(t) = \sum_i \mathbf{1}\{B_i(0) \leq y, B_i(t) > y\} - \sum_i \mathbf{1}\{B_i(0) > y, B_i(t) \leq y\}$$

be the net current of Brownian particles across the point $y$ during time interval $[0, t]$. To have sample paths in $D_{\mathbf{R}}[0, \infty)$, we should replace $Y_\lambda$ with the right-continuous modification $Y_\lambda^+$ defined by $Y_\lambda^+(t) = Y_\lambda(t+)$. This change does not affect finite-dimensional distributions. The covariance of $Y_\lambda$ is

$$EY_\lambda(s)Y_\lambda(t) = \frac{\lambda}{\sqrt{2\pi}}(\sqrt{s} + \sqrt{t} - \sqrt{t-s}) \qquad \text{for } s \leq t.$$

THEOREM 2.4. *As $\lambda \to \infty$, the process $\lambda^{-1/2} Y_\lambda^+(\cdot)$ converges in distribution on the space $D_\mathbf{R}[0,\infty)$ to fractional Brownian motion $Z(\cdot)$ with covariance*

$$EZ(s)Z(t) = \frac{1}{\sqrt{2\pi}}(\sqrt{s} + \sqrt{t} - \sqrt{t-s}) \quad \text{for } s \leq t.$$

The calculations of this paper can be adapted from the random walk situation to the Brownian situation and we omit the explicit proof.

Before turning to the proofs, we want to compare the independent walks with Hammersley's process.

**3. A bound on second-order fluctuations in Hammersley's process.** In this section we look at Hammersley's process from the same perspective from which the previous section studied independent random walks. The difference is that now there is genuine interaction among the particles. The hydrodynamic equation is a nonlinear Hamilton–Jacobi equation. The characteristics of the equation can meet and form shocks. Currently, we cannot prove the equivalent of Theorem 2.2 for Hammersley's process. We can only give bounds on the tails of the second-order fluctuations which suggest that if there is a limit, it should be on the $n^{1/3}$ scale.

The state of Hammersley's process at time $t$ is $z_t = (z_t(i) : i \in \mathbf{Z})$. Depending on the preferred interpretation, variable $z_t(i) \in \mathbf{R}$ is the location of particle labeled $i$ at time $t$ or the height of the interface over site $i$. The dynamics preserves the ordering $z_t(i-1) \leq z_t(i)$. Particles jump to the left, according to this rule. If the state at time $t$ is $z_t = (z_t(i) : i \in \mathbf{Z})$, then particle $i$ has an instantaneous rate $z_t(i) - z_t(i-1)$ of jumping, independently of all other particles. And when particle $i$ jumps, its new location is chosen uniformly at random from the interval $(z_t(i-1), z_t(i))$. This happens independently for all particles $i$.

This process can be defined in terms of a special graphical construction that utilizes the increasing sequences in a space-time Poisson point process, see [1, 15] or [18].

The process of increment variables $\eta_t = (\eta_t(i) : i \in \mathbf{Z})$ is defined as before by

$$\eta_t(i) = z_t(i) - z_t(i-1),$$

and is also known as the "stick process." The dynamics of $\eta_t$ operates as follows. For each $i \in \mathbf{Z}$, at rate equal to $\eta(i)$, this stick-breaking event happens: pick $u$ uniformly distributed on $[0, \eta(i)]$, and replace the state $\eta$ with the new state

$$\eta^{u,i,i+1}(j) = \begin{cases} \eta(i) - u, & j = i, \\ \eta(i+1) + u, & j = i+1, \\ \eta(j), & j \neq i, i+1. \end{cases}$$

In other words, $\eta^{u,i,i+1}$ represents the stick configuration after a piece of size $u$ has been moved from site $i$ to $i+1$. This process can be rigorously defined on a certain subspace of the full product space $[0,\infty)^{\mathbf{Z}}$, see [15] for details.

Next we describe one set of hypotheses under which the hydrodynamic limit and the trivial fluctuations (1.2) can be proved. Then we state a bound on the size of the second-order fluctuations. The setting is again that of a sequence of processes $z_t^n$, $n \in \mathbf{N}$.

ASSUMPTION D. Assume given a nondecreasing Lipschitz function $u_0$ on $\mathbf{R}$ and a bounded, continuous, nonnegative function $v_0$ on $\mathbf{R}$. For all $n$, $z_0^n(0) = 0$ and the initial increment variables $\{\eta_0^n(i) : i \in \mathbf{Z}\}$ are mutually independent with means and variances

$$E[\eta_0^n(i)] = nu_0\left(\frac{i}{n}\right) - nu_0\left(\frac{i-1}{n}\right) \quad \text{and} \quad \operatorname{Var}[\eta_0^n(i)] = v_0\left(\frac{i}{n}\right).$$

Furthermore, assume uniformly bounded fourth moments:

$$\sup_{n \in \mathbf{N}, i \in \mathbf{Z}} E[\eta_0^n(i)^4] < \infty.$$

Let $u(x,t)$ be the unique viscosity solution of the Hamilton–Jacobi equation

(3.1) $$u_t + (u_x)^2 = 0, \qquad u(x,0) = u_0(x).$$

Equivalently, $u$ is defined for $t > 0$ by the Hopf–Lax formula

(3.2) $$u(x,t) = \inf_{y \,:\, y \leq x}\left\{u_0(y) + \frac{(x-y)^2}{4t}\right\}.$$

The hypotheses guarantee that there exists a nonempty compact set $I(x,t) \subseteq (-\infty, x]$ on which the infimum in (3.2) is achieved:

$$I(x,t) = \left\{y \leq x : u(x,t) = u_0(y) + \frac{(x-y)^2}{4t}\right\}.$$

A point $(x,t)$ is a *shock* if $I(x,t)$ is not a singleton. Equivalently, $u(x,t)$ is not differentiable in the $x$ variable at $(x,t)$. For a fixed $t > 0$, there are, at most, countably many shocks. Shocks cannot happen for the linear equation (2.14) of independent particles, because its characteristics are parallel straight lines.

Here is the starting point: the hydrodynamic limit and the fluctuations transported by the characteristics.

THEOREM 3.1. *Assume Assumption* D. *Then, for each* $(x,t) \in \mathbf{R} \times [0,\infty)$,

(3.3) $$\lim_{n \to \infty} n^{-1} z_{nt}^n([nx]) = u(x,t) \qquad \text{in probability}$$

*and*

$$(3.4) \quad \lim_{n\to\infty} \left\{ \frac{z_{nt}^n([nx]) - nu(x,t)}{\sqrt{n}} - \inf_{y \in I(x,t)} \frac{z_0^n([ny]) - nu_0(y)}{\sqrt{n}} \right\} = 0$$

*in probability.*

It is possible to construct the processes $z_t^n$ on the same probability space with a two-sided standard Brownian motion $B(\cdot)$ such that these limits hold in probability, for all $(x,t)$:

$$(3.5) \quad \lim_{n\to\infty} \frac{z_{nt}^n([nx]) - nu(x,t)}{\sqrt{n}} = \inf_{y \in I(x,t)} B\left( \int_0^y v_0(s)\, ds \right).$$

This theorem is proved in [18]. The hypotheses for (3.5) in [18] are more stringent than the ones used above (in [18] the initial increments are assumed exponentially distributed) and the conclusion is stronger (a.s. convergence). The argument in [18] gives convergence in probability in (3.5) under the fourth moment bound included in Assumption D.

The infimum in (3.4) and (3.5) is in some sense the same infimum that appears in the Hopf–Lax formula (3.2), which is inherited by a particle-level variational formulation (5.2).

The result of Baik, Deift and Johansson [2] gives the fluctuations for Hammersley's process from the following particular deterministic initial state: $z_0(i) = 0$ for $i \leq 0$ and $z_0(i) = \infty$ for $i > 0$. In this situation the number of particles in space interval $(0, x]$ at time $t$ equals the maximal number $\mathbf{L}(x,t)$ of space-time Poisson points on an increasing path in the rectangle $(0,x] \times (0,t]$. This connection comes from the graphical construction of Hammersley's process. The distributional limit for $n^{-1/3}\{\mathbf{L}(nx, nt) - 2n\sqrt{xt}\}$ in [2] can be translated into a limit for $n^{-1/3}\{z_{nt}([nx]) - nx^2/(4t)\}$ for $x, t > 0$.

We saw for the independent random walk model that the scale of the fluctuations from deterministic initial conditions is the scale of the second-order fluctuations. So, given the Baik–Deift–Johansson result, we would expect the next order of fluctuations for Hammersley's process at scale $n^{1/3}$. To capture these fluctuations, fix $(x,t) \in \mathbf{R} \times (0,\infty)$ and define

$$Y_n = \{z_{nt}^n([nx]) - nu(x,t)\} - \inf_{y \in I(x,t)} \{z_0^n([ny]) - nu_0(y)\}.$$

We shall prove a bound on the tails of $Y_n$ that suggests $n^{1/3}$ as the correct order. We need one more assumption. Given $(x,t)$, let

$$\Phi(y) = u_0(y) + \frac{(x-y)^2}{4t}$$

be the quantity minimized over $y$ in the Hopf–Lax formula (3.2).

ASSUMPTION E. *Given $(x,t)$, the minimizers in (3.2) are uniformly quadratic, in other words, there exist $c_1, \delta > 0$ such that*

$$\Phi(y) - \Phi(\bar{y}) \geq c_1 (y - \bar{y})^2 \tag{3.6}$$

*for all $y \in \mathbf{R}$ and $\bar{y} \in I(x,t)$ such that $|y - \bar{y}| \leq \delta$.*

THEOREM 3.2. *Let $(x,t) \in \mathbf{R} \times (0, \infty)$. Under Assumptions D and E, the sequence*

$$\left\{ \frac{Y_n}{n^{1/3} \log n} : n \geq 1 \right\}$$

*is tight.*

We turn to proofs, beginning with the random walk model.

## 4. Proofs for the random walk model.

The main work is in proving Theorem 2.2. Along the way we derive an estimate that takes care of Theorem 2.3. Last we explain how Theorem 2.1 follows. We start by proving the convergence of finite-dimensional distributions to the correct limit, and then prove tightness at the process level.

Using (2.20), write $Y_n(t) = Y_{n,1}(t) - Y_{n,2}(t)$, where

$$Y_{n,1}(t) = \sum_i \mathbf{1}\{X_i^n(0) \geq [n\bar{y}] + 1, X_i^n(nt) \leq [n\bar{y}] + [nbt]\} \quad \text{and}$$

$$Y_{n,2}(t) = \sum_i \mathbf{1}\{X_i^n(0) \leq [n\bar{y}], X_i^n(nt) > [n\bar{y}] + [nbt]\}. \tag{4.1}$$

$Y_{n,1}(t)$ and $Y_{n,2}(t)$ represent contributions of slow and fast random walks, respectively. Next write

$$Y_n(t) = (EY_{n,1}(t) - EY_{n,2}(t)) \tag{4.2}$$
$$+ (Y_{n,1}(t) - EY_{n,1}(t)) - (Y_{n,2}(t) - EY_{n,2}(t)).$$

We look at the behavior of these three terms on the scale $n^{1/4}$. Note that $Y_{n,1}(t)$ and $Y_{n,2}(t)$ are independent of each other.

### 4.1. Convergence of finite-dimensional distributions.

PROPOSITION 4.1. *Fix $N$ time points*

$$0 \leq t_1 < t_2 < \cdots < t_N.$$

*As $n \to \infty$ the vector $n^{-1/4}(Y_n(t_1), Y_n(t_2), \ldots, Y_n(t_N))$ converges in distribution to the mean-zero Gaussian random vector $(Z(t_1), Z(t_2), \ldots, Z(t_N))$ with the covariance defined in (2.21).*

PROOF. Since $Y_n(0)$ is identically zero, we may as well assume that $t_1 > 0$. By the Cramér–Wold device, it suffices to show the convergence of the linear combination

$$(4.3) \qquad n^{-1/4} \sum_{i=1}^{N} \theta_i Y_n(t_i)$$

for an arbitrary vector $\theta = (\theta_1, \ldots, \theta_N) \in \mathbf{R}^N$ and arbitrary $N \in \mathbf{N}$.

In Lemma 4.6 below we show that $n^{-1/4} E Y_n(t) = n^{-1/4}(EY_{n,1}(t) - EY_{n,2}(t))$ vanishes as $n \to \infty$. Using the decomposition (4.2) and ignoring the first term that vanishes, what we actually prove is the weak convergence of the difference

$$(4.4) \; n^{-1/4} \sum_{i=1}^{N} \theta_i(Y_{n,1}(t_i) - EY_{n,1}(t_i)) - n^{-1/4} \sum_{i=1}^{N} \theta_i(Y_{n,2}(t_i) - EY_{n,2}(t_i)).$$

Since the two sums above are independent, we can treat them separately. We shall show below that they converge to mean-zero normal distributions with variances

$$\sigma_1^2 = \sum_{i=1}^{N} \theta_i^2 \sqrt{\kappa_2} \bigg\{ \rho_0(\bar{y}) \int_{-\infty}^{0} P(B_{t_i} > z) P(B_{t_i} \leq z) \, dz$$
$$+ v_0(\bar{y}) \int_{-\infty}^{0} P(B_{t_i} \leq z)^2 \, dz \bigg\}$$
$$(4.5) \qquad + 2 \sum_{1 \leq i < j \leq N} \theta_i \theta_j \sqrt{\kappa_2} \bigg\{ \rho_0(\bar{y}) \int_{-\infty}^{0} [P(B_{t_i} > z) P(B_{t_j} \leq z)$$
$$- P(B_{t_i} > z \geq B_{t_j})] \, dz$$
$$+ v_0(\bar{y}) \int_{-\infty}^{0} P(B_{t_i} \leq z) P(B_{t_j} \leq z) \, dz \bigg\}$$

and

$$\sigma_2^2 = \sum_{i=1}^{N} \theta_i^2 \sqrt{\kappa_2} \bigg\{ \rho_0(\bar{y}) \int_{0}^{\infty} P(B_{t_i} > z) P(B_{t_i} \leq z) \, dz$$
$$+ v_0(\bar{y}) \int_{0}^{\infty} P(B_{t_i} > z)^2 \, dz \bigg\}$$
$$(4.6) \qquad + 2 \sum_{1 \leq i < j \leq N} \theta_i \theta_j \sqrt{\kappa_2} \bigg\{ \rho_0(\bar{y}) \int_{0}^{\infty} [P(B_{t_i} > z) P(B_{t_j} \leq z)$$
$$- P(B_{t_i} > z \geq B_{t_j})] \, dz$$
$$+ v_0(\bar{y}) \int_{0}^{\infty} P(B_{t_i} > z) P(B_{t_j} > z) \, dz \bigg\}.$$

$B_t$ above stands for standard one-dimensional Brownian motion. The quantities $\sigma_1^2$ and $\sigma_2^2$ are actually equal. Their different formulas connect to the calculations by which they arise in the proof. The limit of (4.3) is then a centered Gaussian with variance $\sigma_1^2 + \sigma_2^2$. To evaluate this sum, note that the terms with $\rho_0(\bar{y})$ add up to integrals over $(-\infty, \infty)$, while the terms with $v_0(\bar{y})$ are actually equal. We leave the proof of the next lemma to the reader (calculus and Fubini's theorem are needed).

LEMMA 4.1. *For $0 \leq s \leq t$, we have these formulas:*

$$\text{(4.7)} \qquad \int_{-\infty}^{\infty} P(B_s > z) P(B_t \leq z) \, dz = \frac{\sqrt{s+t}}{\sqrt{2\pi}},$$

$$\text{(4.8)} \qquad \int_0^{\infty} P(B_s > z) P(B_t > z) \, dz = \frac{\sqrt{s} + \sqrt{t} - \sqrt{s+t}}{2\sqrt{2\pi}}$$

*and*

$$\text{(4.9)} \qquad \int_{-\infty}^{\infty} P(B_s > z \geq B_t) \, dz = \frac{\sqrt{t-s}}{\sqrt{2\pi}}.$$

From this lemma and the definitions (4.5)–(4.6), we get

$$\sigma_1^2 + \sigma_2^2 = \sum_{i=1}^{N} \theta_i^2 \sqrt{\kappa_2} \left\{ \rho_0(\bar{y}) \frac{\sqrt{t_i}}{\sqrt{\pi}} + v_0(\bar{y}) \frac{(\sqrt{2}-1)\sqrt{t_i}}{\sqrt{\pi}} \right\}$$
$$+ 2 \sum_{1 \leq i < j \leq N} \theta_i \theta_j \sqrt{\kappa_2} \left\{ \rho_0(\bar{y}) \left( \frac{\sqrt{t_i + t_j}}{\sqrt{2\pi}} - \frac{\sqrt{t_j - t_i}}{\sqrt{2\pi}} \right) \right.$$
$$\left. + v_0(\bar{y}) \frac{\sqrt{t_i} + \sqrt{t_j} - \sqrt{t_i + t_j}}{\sqrt{2\pi}} \right\}$$
$$= \sum_{1 \leq i, j \leq N} \theta_i \theta_j EZ(t_i) Z(t_j),$$

where the last equality comes from (2.21). Thus, the linear combination in (4.3) converges in distribution to the linear combination $\sum_{i=1}^{N} \theta_i Z(t_i)$. Since the vector $\theta$ was arbitrary, Proposition 4.1 follows.

Equations (4.7)–(4.9) can be manipulated to show that the functions $\sqrt{s+t} - \sqrt{s \vee t - s \wedge t}$ and $\sqrt{s} + \sqrt{t} - \sqrt{s+t}$ are positive definite. This ensures that (2.21) is a legitimate covariance of a Gaussian process for all nonnegative values $\rho_0(\bar{y})$ and $v_0(\bar{y})$.

It remains now to prove the weak convergence of the sums in (4.4) and the vanishing of the mean $n^{-1/4} EY_n(t)$ in the limit.

LEMMA 4.2. *As $n \to \infty$, $n^{-1/4} \sum_{i=1}^{N} \theta_i (Y_{n,1}(t_i) - EY_{n,1}(t_i))$ converges weakly to a mean-zero normal distribution with variance $\sigma_1^2$ defined by (4.5).*

Lemma 4.2 will be proved after some preliminary steps. Relabel the random walks so that $X_{m,j}(\cdot)$ is the $j$th random walk that starts at site $[n\bar{y}]+m$. Then

$$Y_{n,1}(t) = \sum_{m=1}^{\infty} \sum_{j=1}^{\eta_0^n([n\bar{y}]+m)} \mathbf{1}\{X_{m,j}(nt) - X_{m,j}(0) \leq [nbt] - m\}.$$

Since the random walks are independent of the initial occupation numbers $\eta_0^n(x)$,

(4.10) $$EY_{n,1}(t) = \sum_{m=1}^{\infty} \rho_0^n([n\bar{y}]+m) \cdot P\{X(nt) \leq [nbt] - m\},$$

where $X(\cdot)$ represents a random walk with rates $p(x)$ starting at the origin.

Write

$$n^{-1/4}(Y_{n,1}(t) - EY_{n,1}(t)) = \sum_{m=1}^{\infty} U_m(t),$$

with mean zero summands

$$U_m(t) = n^{-1/4} \sum_{j=1}^{\eta_0^n([n\bar{y}]+m)} \mathbf{1}\{X_{m,j}(nt) - X_{m,j}(0) \leq [nbt] - m\}$$
$$- n^{-1/4} \rho_0^n([n\bar{y}]+m) \cdot P\{X(nt) \leq [nbt] - m\}.$$

For fixed $n$ and $t$, the variables $\{U_m(t) : m \geq 1\}$ are independent. Abbreviate

$$\overline{U}_m = \sum_{i=1}^{N} \theta_i U_m(t_i).$$

Rearrange as follows:

$$n^{-1/4} \sum_{i=1}^{N} \theta_i (Y_{n,1}(t_i) - EY_{n,1}(t_i)) = \sum_{m=1}^{\infty} \sum_{i=1}^{N} \theta_i U_m(t_i) = \sum_{m=1}^{\infty} \overline{U}_m = S_1 + S_2,$$

where

$$S_1 = \sum_{m=1}^{[r(n)\sqrt{n}]} \overline{U}_m \quad \text{and} \quad S_2 = \sum_{m=[r(n)\sqrt{n}]+1}^{\infty} \overline{U}_m.$$

We shall apply the Lindeberg–Feller theorem to $S_1$ and show that $S_2 \to 0$ in $L^2$. To this end, we make $r(n) \nearrow \infty$ sufficiently slowly. Let

$$H(M) = \sup_{n \geq 1, x \in \mathbf{Z}} E[\eta_0^n(x)^2 \mathbf{1}\{\eta_0^n(x) \geq M\}].$$

$H(M) \to 0$ as $M \to \infty$ by assumption (2.8). Now choose $r(n) = o(\sqrt{\log n})$ so that $r(n) \nearrow \infty$ while

$$r(n)H(n^{1/8}) \to 0. \tag{4.11}$$

The condition $r(n) = o(\sqrt{\log n})$ is imposed so that later we can use assumption (2.9).

We first show $ES_2^2 \to 0$. $S_2$ is a sum of independent mean-zero terms $\overline{U}_m$, and so

$$ES_2^2 = \sum_{m=[r(n)\sqrt{n}]+1}^{\infty} E\overline{U}_m^2 \leq \|\theta\|^2 \sum_{i=1}^{N} \sum_{m=[r(n)\sqrt{n}]+1}^{\infty} EU_m^2(t_i).$$

We wrote $\|\theta\|$ for the Euclidean norm and used the Schwarz inequality. Since $N$ is fixed it suffices to show that for a fixed $t$,

$$\lim_{n \to \infty} \sum_{m=[r(n)\sqrt{n}]+1}^{\infty} EU_m^2(t) = 0. \tag{4.12}$$

Recall that the variance of a random sum $T_K = \sum_{i=1}^{K} Z_i$ with i.i.d. summands $Z_i$ independent of the random $K$ is

$$\mathrm{Var}[T_K] = EK \cdot \mathrm{Var}\, Z_1 + (EZ_1)^2 \cdot \mathrm{Var}\, K. \tag{4.13}$$

The variance of the indicator $\mathbf{1}\{X(nt) \leq [nbt] - m\}$ is $P\{X(nt) \leq [nbt] - m\}P\{X(nt) > [nbt] - m\}$. $X(nt)$ is a sum of a Poisson($nt$) distributed number of independent jumps, each jump distributed according to $\{p(x)\}$. Hence, the variance of $X(nt)$ is $nt\kappa_2$. By Donsker's invariance principle, the process $\{(X(nt) - [nbt])/\sqrt{n\kappa_2} : t \geq 0\}$ converges weakly to standard one-dimensional Brownian motion.

By the definition of $U_m(t)$ and (4.13),

$$E[U_m^2(t)] = \mathrm{Var}\left[n^{-1/4} \sum_{j=1}^{\eta_0^n([n\bar{y}]+m)} \mathbf{1}\{X_{m,j}(nt) - X_{m,j}(0) \leq [nbt] - m\}\right]$$

$$\begin{aligned}(4.14) \quad &= n^{-1/2}\rho_0^n([n\bar{y}]+m)P\{X(nt) \leq [nbt]-m\}P\{X(nt) > [nbt]-m\} \\ &\quad + n^{-1/2}v_0^n([n\bar{y}]+m)P\{X(nt) \leq [nbt]-m\}^2.\end{aligned}$$

In particular, we get the bound

$$E[U_m^2(t)] \leq Cn^{-1/2}P\{X(nt) \leq [nbt] - m\} \tag{4.15}$$

by the uniform bound (2.8) on the moments.

By standard large deviation theory and assumption (2.6), as $s \to \infty$, the random walk $X(s)$ has a rate function $I$ which is convex and quadratic around its unique minimum at $b$. Consequently, for arbitrarily small $\alpha > 0$, there exists a constant $0 < K < \infty$ such that

$$(4.16) \qquad I(b+z) \geq \begin{cases} Kz^2, & |z| \leq \alpha, \\ K|z|, & |z| \geq \alpha. \end{cases}$$

The rate function gives the bounds

$$P\{X(s) \leq sb - su\} \leq \exp\{-sI(b-u)\}$$

and

$$P\{X(s) \geq sb + su\} \leq \exp\{-sI(b+u)\}$$

for all $u \geq 0$ and $s > 0$. Even though a large deviations rate function is an asymptotic notion, these bounds are valid already for finite $s$ by virtue of superadditivity. We use these first in the form

$$(4.17) \qquad \begin{aligned} P\{X(nt) \leq [nbt] - m\} &\leq \exp\left\{-ntI\left(b - \frac{m}{nt}\right)\right\} \\ &\leq \begin{cases} \exp\left\{-K\dfrac{m^2}{nt}\right\}, & 0 \leq m \leq nt\alpha, \\ \exp\{-Km\}, & m \geq nt\alpha. \end{cases} \end{aligned}$$

Consequently, applying (4.15) and (4.17),

$$\sum_{m=[r(n)\sqrt{n}]+1}^{\infty} E[U_m(t)^2]$$

$$\leq Cn^{-1/2} \sum_{m=[r(n)\sqrt{n}]+1}^{[nt\alpha]} \exp\{-Km^2(nt)^{-1}\} + Cn^{-1/2} \sum_{m=[nt\alpha]+1}^{\infty} \exp\{-Km\}$$

$$\leq C \int_{r(n)}^{\infty} e^{-Kt^{-1}x^2}\, dx + \frac{C}{1 - e^{-K}} \cdot n^{-1/2} e^{-K\alpha nt}.$$

This vanishes as $n \to \infty$ due to $r(n) \to \infty$. We have proved (4.12).

Next follows the application of Lindeberg–Feller to $S_1$. Two things need to be checked, namely, that

$$(4.18) \qquad \lim_{n\to\infty} \sum_{m=1}^{[r(n)\sqrt{n}]} E[\overline{U}_m^2 \mathbf{1}\{|\overline{U}_m| \geq \varepsilon\}] = 0$$

for any $\varepsilon > 0$, and, second, that

$$(4.19) \qquad \lim_{n\to\infty} \sum_{m=1}^{[r(n)\sqrt{n}]} E[\overline{U}_m^2] = \sigma_1^2.$$

(See Theorem (4.5) on page 116 in [5].) Note that

$$|U_m(t)| \leq n^{-1/4}\eta_0^n([n\bar{y}] + m) + n^{-1/4}\rho_0^n([n\bar{y}] + m)$$

and so by the uniform bound on moments,

$$|\overline{U}_m| \leq Cn^{-1/4}(\eta_0^n([n\bar{y}] + m) + 1).$$

(4.18) follows from this and property (4.11) of $r(n)$. In (4.11) we used $n^{1/8}$ simply because for any $\varepsilon > 0$, $\varepsilon n^{1/4} > n^{1/8}$ for large enough $n$.

We turn to verify (4.19).

$$\overline{U}_m^2 = \sum_{i=1}^{N} \theta_i^2 U_m^2(t_i) + 2 \sum_{1 \leq i < j \leq N} \theta_i \theta_j U_m(t_i) U_m(t_j),$$

and so the sum in (4.19) can be expressed as

$$(4.20) \quad \sum_{m=1}^{[r(n)\sqrt{n}]} E[\overline{U}_m^2] = \sum_{i=1}^{N} \theta_i^2 S_{1,1}(t_i) + 2 \sum_{1 \leq i < j \leq N} \theta_i \theta_j S_{1,2}(t_i, t_j),$$

where we abbreviated

$$S_{1,1}(t) = \sum_{m=1}^{[r(n)\sqrt{n}]} E[U_m^2(t)]$$

and

$$S_{1,2}(s,t) = \sum_{m=1}^{[r(n)\sqrt{n}]} E[U_m(s)U_m(t)].$$

LEMMA 4.3. *For $t > 0$,*

$$\lim_{n \to \infty} S_{1,1}(t) = \rho_0(\bar{y})\sqrt{\kappa_2} \int_{-\infty}^{0} P(B_t \leq z) P(B_t > z) \, dz$$
$$+ v_0(\bar{y})\sqrt{\kappa_2} \int_{-\infty}^{0} P(B_t \leq z)^2 \, dz,$$

*where $B_t$ is standard one-dimensional Brownian motion.*

PROOF. By (4.14),

$$\sum_{m=1}^{[r(n)\sqrt{n}]} E[U_m^2(t)] = n^{-1/2} \sum_{m=1}^{[r(n)\sqrt{n}]} \rho_0^n([n\bar{y}] + m) P\{X(nt) \leq [nbt] - m\}$$

$$\times P\{X(nt) > [nbt] - m\}$$

(4.21)
$$+ n^{-1/2} \sum_{m=1}^{[r(n)\sqrt{n}]} v_0^n([n\bar{y}] + m) P\{X(nt) \leq [nbt] - m\}^2$$
$$\equiv T_{1,1} + T_{1,2}.$$

The last equality above defines the sums $T_{1,1}$ and $T_{1,2}$. We work with $T_{1,1}$, and leave the analogous arguments for $T_{1,2}$ to the reader.

We bound $T_{1,1}$ from below. In the calculation that follows, $L = L(n)$ is the integer that appeared in assumption (2.9). The $o(n^{-1/4})$ error term below that comes from that assumption is uniform over $k$ because $kL = O(r(n)\sqrt{n}) = O(\sqrt{n\log n})$, which is permitted in assumption (2.9),

$$T_{1,1} = n^{-1/2} \sum_{m=1}^{[r(n)\sqrt{n}]} \rho_0^n([n\bar{y}] + m) P\{X(nt) \leq [nbt] - m\} P\{X(nt) > [nbt] - m\}$$

$$\geq n^{-1/2} \sum_{k=0}^{[L^{-1}r(n)\sqrt{n}]-1} \sum_{j=1}^{L} \rho_0^n([n\bar{y}] + kL + j)$$
$$\times P\{X(nt) \leq [nbt] - (k+1)L\} P\{X(nt) > [nbt] - kL\}$$

$$\geq n^{-1/2} \sum_{k=0}^{[L^{-1}r(n)\sqrt{n}]-1} (L\rho_0(\bar{y}) + L \cdot o(n^{-1/4}))$$
$$\times P\{X(nt) \leq [nbt] - (k+1)L\} P\{X(nt) > [nbt] - kL\}$$

$$\geq \rho_0(\bar{y}) n^{-1/2} L \sum_{k=0}^{[L^{-1}r(n)\sqrt{n}]-1} P\{X(nt) \leq [nbt] - (k+1)L\}$$
$$\times P\{X(nt) > [nbt] - kL\} + o(n^{-1/4}) \cdot O(r(n))$$

$$= \rho_0(\bar{y}) n^{-1/2} L \sum_{k=0}^{[L^{-1}r(n)\sqrt{n}]-1} P\left\{\frac{X(nt) - [nbt]}{\sqrt{n\kappa_2}} \leq -(k+1)\frac{L}{\sqrt{n\kappa_2}}\right\}$$
$$\times P\left\{\frac{X(nt) - [nbt]}{\sqrt{n\kappa_2}} > -k\frac{L}{\sqrt{n\kappa_2}}\right\}$$
$$+ o(n^{-1/4}) \cdot O(r(n)).$$

The last term $o(n^{-1/4}) \cdot O(r(n)) \to 0$ as $n \to \infty$ because $r(n) = o(\sqrt{\log n})$. As $n \to \infty$, a Riemann sum argument, together with the large deviation bounds (4.17), shows that the main part of the lower bound converges to

$$\rho_0(\bar{y}) \sqrt{\kappa_2} \int_{-\infty}^{0} P(B_t \leq z) P(B_t > z) \, dz.$$

Similarly, one derives an upper bound for $T_{1,1}$ with the same limit. This proves the convergence of $T_{1,1}$.

We leave the similar argument for $T_{1,2}$ to the reader. This completes the proof of Lemma 4.3. □

LEMMA 4.4. *For $0 < s < t$,*

$$\lim_{n \to \infty} S_{1,2}(s,t)$$

$$(4.22) \quad = \rho_0(\bar{y})\sqrt{\kappa_2} \int_{-\infty}^{0} \{P(B_s > z)P(B_t \leq z) - P(B_s > z \geq B_t)\}\, dz$$

$$+ v_0(\bar{y})\sqrt{\kappa_2} \int_{-\infty}^{0} P(B_s \leq z)P(B_t \leq z)\, dz.$$

PROOF. Formula (4.13) generalizes in the following way. Assume the i.i.d. random variables $\{Z_i\}$ are independent of the random nonnegative integer $K$, and $f$ and $g$ are bounded measurable functions on the state space of the $\{Z_i\}$. Then

$$\operatorname{Cov}\left[\sum_{i=1}^{K} f(Z_i), \sum_{i=1}^{K} g(Z_i)\right]$$

$$(4.23) \quad = E\left[\left\{\sum_{i=1}^{K} f(Z_i) - EK \cdot Ef(Z_1)\right\}\left\{\sum_{i=1}^{K} g(Z_i) - EK \cdot Eg(Z_1)\right\}\right]$$

$$= EK \cdot \operatorname{Cov}[f(Z_1), g(Z_1)] + \operatorname{Var} K \cdot Ef(Z_1) \cdot Eg(Z_1).$$

Applying this gives

$$EU_m(s)U_m(t)$$
$$= n^{-1/2}\rho_0^n([n\bar{y}] + m)(P\{X(ns) \leq [nbs] - m, X(nt) \leq [nbt] - m\}$$
$$\qquad - P\{X(ns) \leq [nbs] - m\} \cdot P\{X(nt) \leq [nbt] - m\})$$
$$+ n^{-1/2}v_0^n([n\bar{y}] + m)P\{X(ns) \leq [nbs] - m\} \cdot P\{X(nt) \leq [nbt] - m\}.$$

The probabilities in the first term can be rearranged as follows:

$$P\{X(ns) \leq [nbs] - m, X(nt) \leq [nbt] - m\}$$
$$- P\{X(ns) \leq [nbs] - m\} \cdot P\{X(nt) \leq [nbt] - m\}$$
$$= P\{X(ns) > [nbs] - m\} \cdot P\{X(nt) \leq [nbt] - m\}$$
$$- P\{X(ns) > [nbs] - m, X(nt) \leq [nbt] - m\}.$$

With $m = [z\sqrt{n\kappa_2}]$ this last expression converges to

$$P(B_s > z)P(B_t \leq z) - P(B_s > z \geq B_t),$$

which is the integrand of the first integral in (4.22). This points the way, and one can follow the reasoning of the proof of Lemma 4.3. □

Together with (4.20), Lemmas 4.3 and 4.4 prove the limit in (4.19), and, thereby, the central limit theorem claimed in Lemma 4.2. Next we need the corresponding result for the second sum in the difference (4.4).

LEMMA 4.5. $n^{-1/4} \sum_{i=1}^{N} \theta_i (Y_{n,2}(t_i) - EY_{n,2}(t_i))$ *converges weakly to a mean-zero normal distribution with variance $\sigma_2^2$ defined in* (4.6).

PROOF. We have the same argument in principle,

$$Y_{n,2}(t) = \sum_{m=0}^{\infty} \sum_{j=1}^{\eta_0^n([n\bar{y}]-m)} \mathbf{1}\{X_{m,j}(nt) - X_{m,j}(0) > [nbt] + m\}.$$

Then

$$(4.24) \qquad EY_{n,2}(t) = \sum_{m=0}^{\infty} \rho_0^n([n\bar{y}] - m) \cdot P\{X(nt) > [nbt] + m\}.$$

Next write

$$n^{-1/4}(Y_{n,2}(t) - EY_{n,2}(t)) = \sum_{m=0}^{\infty} V_m(t),$$

with independent, mean zero summands

$$V_m(t) = n^{-1/4} \sum_{j=1}^{\eta_0^n([n\bar{y}]-m)} \mathbf{1}\{X_{m,j}(nt) - X_{m,j}(0) > [nbt] + m\}$$
$$- n^{-1/4} \rho_0^n([n\bar{y}] - m) \cdot P\{X(nt) > [nbt] + m\}.$$

With

$$\overline{V}_m = \sum_{i=1}^{N} \theta_i V_m(t_i),$$

first separate out the part

$$S_2 = \sum_{m=[r(n)\sqrt{n}]}^{\infty} \overline{V}_m.$$

Use large deviation estimates to show that $S_2 \to 0$ in $L^2$ as $n \to \infty$. To the remaining part

$$S_1 = \sum_{m=0}^{[r(n)\sqrt{n}]-1} \overline{V}_m,$$

apply Lindeberg–Feller. The details are similar to those in the proof of Lemma 4.2. □

LEMMA 4.6. *For any $0 < T < \infty$,*

$$\lim_{n\to\infty} \sup_{0\le t\le T} n^{-1/4}|EY_n(t)| = \lim_{n\to\infty} \sup_{0\le t\le T} n^{-1/4}|EY_{n,1}(t) - EY_{n,2}(t)| = 0.$$

PROOF. From (4.10) and (4.24),

$$EY_{n,1}(t) - EY_{n,2}(t) = S + R_1 - R_2,$$

where

$$S = \rho_0(\bar{y})\left\{\sum_{m=1}^{\infty} P\{X(nt) \le [nbt] - m\} - \sum_{m=0}^{\infty} P\{X(nt) > [nbt] + m\}\right\},$$

$$R_1 = \sum_{m=1}^{\infty} (\rho_0^n([n\bar{y}] + m) - \rho_0(\bar{y}))P\{X(nt) \le [nbt] - m\}$$

and

$$R_2 = \sum_{m=0}^{\infty} (\rho_0^n([n\bar{y}] - m) - \rho_0(\bar{y}))P\{X(nt) > [nbt] + m\}.$$

The part of $S$ in braces equals $[nbt] - EX(nt) = [nbt] - nbt$. Thus, $|S| \le \rho_0(\bar{y})$ uniformly over $n$ and $t$.

Next we show $R_1 = o(n^{1/4})$ uniformly over $t \in [0, T]$. $R_1$ is a sum of bounded terms, and any sum of bounded terms can be rearranged in this manner:

$$\sum_{m=1}^{\infty} a_m = \frac{1}{L}\sum_{j=1}^{L}\sum_{m=1}^{\infty} a_m$$

$$= \frac{1}{L}\sum_{j=1}^{L}\sum_{m=1}^{j-1} a_m + \frac{1}{L}\sum_{j=1}^{L}\sum_{m=j}^{M+j-1} a_m + \frac{1}{L}\sum_{j=1}^{L}\sum_{m=M+j}^{\infty} a_m$$

$$= O(L) + \sum_{m=1}^{M}\frac{1}{L}\sum_{j=0}^{L-1} a_{m+j} + \frac{1}{L}\sum_{j=1}^{L}\sum_{m=M+j}^{\infty} a_m,$$

and so

$$\left|\sum_{m=1}^{\infty} a_m - \sum_{m=1}^{M}\frac{1}{L}\sum_{j=0}^{L-1} a_{m+j}\right| \le O(L) + \sum_{m=M+1}^{\infty} |a_m|.$$

In our situation $L = o(n^{1/4})$ from assumption (2.9). We take $M = M(n) = [c\sqrt{n \log n}\,]$ for a large enough constant $c$. Then the large deviation estimates (4.17) show that

$$\sup_{0 \le t \le T} \sum_{m > c\sqrt{n \log n}} P\{X(nt) \le [nbt] - m\} \to 0 \qquad \text{as } n \to \infty.$$

It remains to show that the sum

$$R_{1,1} = \sum_{m=1}^{M} \frac{1}{L} \sum_{j=0}^{L-1} (\rho_0^n([n\bar{y}] + m + j) - \rho_0(\bar{y})) P\{X(nt) \le [nbt] - m - j\}$$

is $o(n^{1/4})$ uniformly over $t$. Rewrite $R_{1,1}$ as

$$R_{1,1} = \sum_{m=1}^{M} \frac{1}{L} \sum_{j=0}^{L-1} (\rho_0^n([n\bar{y}] + m + j) - \rho_0(\bar{y})) \sum_{k=j}^{\infty} P\{X(nt) = [nbt] - m - k\}$$

$$= \sum_{m=1}^{M} \sum_{k=0}^{\infty} P\{X(nt) = [nbt] - m - k\} \frac{1}{L} \sum_{j=0}^{(L-1)\wedge k} (\rho_0^n([n\bar{y}] + m + j) - \rho_0(\bar{y}))$$

$$= \sum_{m=1}^{M} \sum_{k=0}^{L-1} P\{X(nt) = [nbt] - m - k\} \frac{1}{L} \sum_{j=0}^{k} (\rho_0^n([n\bar{y}] + m + j) - \rho_0(\bar{y}))$$

$$+ \sum_{m=1}^{M} P\{X(nt) \le [nbt] - m - L\} \frac{1}{L} \sum_{j=0}^{L-1} (\rho_0^n([n\bar{y}] + m + j) - \rho_0(\bar{y})).$$

The next to last line above is $O(L) = o(n^{1/4})$, as can be seen by replacing $\rho_0^n([n\bar{y}] + m + j) - \rho_0(\bar{y})$ with a uniform upper bound and then summing the probabilities over $m$. For the last line, use assumption (2.9) to replace each $L^{-1} \sum_{j=0}^{L-1}(\rho_0^n([n\bar{y}] + m + j) - \rho_0(\bar{y}))$ with $o(n^{-1/4})$ uniformly over $m$. Here the assumption $m = O(\sqrt{n \log n}\,)$ is used. Then we have

$$o(n^{-1/4}) \cdot \sum_{m=1}^{M} P\{X(nt) \le [nbt] - m - L\} = o(n^{-1/4}) \cdot O(n^{1/2}) = o(n^{1/4}).$$

The last $O(n^{1/2})$ bound comes from

$$\sum_{m=1}^{\infty} P\{X(nt) \le [nbt] - m\} = E[([nbt] - X(nt))_+] \le C_1 + C_2 n^{1/2}$$

uniformly over $t \ge 0$ and $n$. This is a consequence of the central limit theorem and uniform integrability from assumption (2.6).

We leave the similar treatment of $R_2$ to the reader. $\square$

Proposition 4.1 is now proved, and we turn to tightness at the process level. $\square$

4.2. *Moment bound for time increment.* Using representation (2.20), the difference $Y_n(t) - Y_n(s)$ for $0 < s < t$ simplifies to

$$Y_n(t) - Y_n(s) = \sum_i (\mathbf{1}\{X_i^n(ns) > [n\bar{y}] + [nbs], X_i^n(nt) \leq [n\bar{y}] + [nbt]\}$$

$$- \mathbf{1}\{X_i^n(ns) \leq [n\bar{y}] + [nbs], X_i^n(nt) > [n\bar{y}] + [nbt]\}).$$

Let

$$\overline{Y}_n(t) = Y_n(t) - EY_n(t).$$

By virtue of Lemma 4.6, weak convergence of the processes $n^{-1/4}\overline{Y}_n(\cdot)$ is equivalent to weak convergence of $n^{-1/4}Y_n(\cdot)$. We shall, in fact, work with the centered processes $\overline{Y}_n(\cdot)$.

Relabel the random walks so that $X_{m,j}^n(\cdot)$, $1 \leq j \leq \eta_0^n(m)$, are the random walks that start at site $m$ in the process indexed by $n$. Let

$$A_{m,j} = \{X_{m,j}^n(ns) > [n\bar{y}] + [nbs], X_{m,j}^n(nt) \leq [n\bar{y}] + [nbt]\},$$
$$B_{m,j} = \{X_{m,j}^n(ns) \leq [n\bar{y}] + [nbs], X_{m,j}^n(nt) > [n\bar{y}] + [nbt]\}$$

and

$$G_m = \sum_{j=1}^{\eta_0^n(m)} (\mathbf{1}_{A_{m,j}} - \mathbf{1}_{B_{m,j}}) - \rho_0^n(m)(P(A_{m,1}) - P(B_{m,1})).$$

Then

$$\overline{Y}_n(t) - \overline{Y}_n(s) = \sum_{m \in \mathbf{Z}} G_m$$

is a sum of independent mean-zero random variables. Recall now assumption (2.8), according to which $E[\eta_0^n(x)^6]$ is uniformly bounded over $n$ and $x$.

PROPOSITION 4.2. *There exists a constant $C$ such that for all $n \in \mathbf{N}$ and $0 \leq s < t$,*

(4.25) $\quad E[(\overline{Y}_n(t) - \overline{Y}_n(s))^6] \leq C(n^{1/2}(t-s)^{1/2} + n^{3/2}(t-s)^{3/2} + 1).$

PROOF. To prove Proposition 4.2, start with

$$E[(\overline{Y}_n(t) - \overline{Y}_n(s))^6]$$

(4.26)
$$= \sum_{m_1,\ldots,m_6 \in \mathbf{Z}} E[G_{m_1} G_{m_2} \cdots G_{m_6}]$$

$$= \sum_m E[G_m^6] + \binom{6}{4} \sum_{m_1 \neq m_2} E[G_{m_1}^4] E[G_{m_2}^2]$$

$$+ \binom{6}{3} \sum_{m_1 < m_2} E[G_{m_1}^3] E[G_{m_2}^3]$$

$$+ \binom{6}{2\ 2\ 2} \sum_{m_1 < m_2 < m_3} E[G_{m_1}^2] E[G_{m_2}^2] E[G_{m_3}^2]$$

$$\leq C \Bigg\{ \sum_{m \in \mathbf{Z}} E[G_m^6] + \bigg( \sum_{m \in \mathbf{Z}} EG_m^4 \bigg) \bigg( \sum_{m \in \mathbf{Z}} EG_m^2 \bigg)$$

$$+ \bigg( \sum_{m \in \mathbf{Z}} |EG_m^3| \bigg)^2 + \bigg( \sum_{m \in \mathbf{Z}} EG_m^2 \bigg)^3 \Bigg\}.$$

Above we collected terms $E[G_{m_1} G_{m_2} \cdots G_{m_6}]$ according to how many times distinct sites appear among the indices $m_1, \ldots, m_6$. Independence and $EG_m = 0$ eliminate all terms $E[G_{m_1} G_{m_2} \cdots G_{m_6}]$ where an index appears by itself. This point is actually critical for obtaining (4.25).

LEMMA 4.7. *There exists a constant $C$ such that for each positive integer $1 \leq k \leq 6$ and for all $m$,*

$$E[|G_m|^k] \leq C(P(A_{m,1}) + P(B_{m,1})).$$

PROOF.

$$E[|G_m|^k] \leq E\Bigg[\bigg(\sum_{j=1}^{\eta_0^n(m)} |\mathbf{1}_{A_{m,j}} - \mathbf{1}_{B_{m,j}}| + \rho_0^n(m)|P(A_{m,1}) - P(B_{m,1})|\bigg)^k\Bigg]$$

$$= \sum_{i=0}^{k} \binom{k}{i} E\Bigg[\bigg(\sum_{j=1}^{\eta_0^n(m)} |\mathbf{1}_{A_{m,j}} - \mathbf{1}_{B_{m,j}}|\bigg)^i\Bigg]$$

$$\times \rho_0^n(m)^{k-i} |P(A_{m,1}) - P(B_{m,1})|^{k-i}.$$

The terms with $i < k$ are bounded by $C(P(A_{m,1}) + P(B_{m,1}))$ as the conclusion demands. For the $i = k$ term we bound as follows: since $|P(A_{m,j}) - P(B_{m,j})| \leq 1$,

$$E\Bigg[\bigg(\sum_{j=1}^{\eta_0^n(m)} |\mathbf{1}_{A_{m,j}} - \mathbf{1}_{B_{m,j}}|\bigg)^k\Bigg]$$

$$= \sum_{h=1}^{\infty} P\{\eta_0^n(m) = h\} \cdot E\Bigg[\bigg(\sum_{j=1}^{h} |\mathbf{1}_{A_{m,j}} - \mathbf{1}_{B_{m,j}}|\bigg)^k\Bigg]$$

$$\leq \sum_{h=1}^{\infty} P\{\eta_0^n(m) = h\} h^k E[|\mathbf{1}_{A_{m,1}} - \mathbf{1}_{B_{m,1}}|]$$

$$\leq E[\eta_0^n(m)^k](P(A_{m,1}) + P(B_{m,1})). \qquad \square$$

We can now bound all the sums on line (4.26) by

$$(4.27) \qquad \sum_{m\in\mathbf{Z}} E|G_m|^k \leq C \sum_{m\in\mathbf{Z}} P(A_{m,1}) + C \sum_{m\in\mathbf{Z}} P(B_{m,1}),$$

so next we estimate the sums of probabilities on the right.

LEMMA 4.8. *There is a finite constant $C$ such that, for all $0 \leq s \leq t$,*

$$(4.28) \qquad \sum_{m\in\mathbf{Z}} E|G_m|^k \leq C(\sqrt{n(t-s)} + 1).$$

PROOF. Write $X(\cdot)$ for a representative random walk that starts at the origin.

$$\sum_{m\in\mathbf{Z}} P(A_{m,1}) = \sum_{m\in\mathbf{Z}} P\{X_{m,1}^n(ns) > [n\bar{y}] + [nbs], X_{m,1}^n(nt) \leq [n\bar{y}] + [nbt]\}$$

$$= \sum_{m\in\mathbf{Z}} P\{X(ns) > [nbs] - m, X(nt) \leq [nbt] - m\}.$$

Apply the Markov property to turn the sum into

$$\sum_{m\in\mathbf{Z}} \sum_{\ell<m} P\{X(ns) = [nbs] - \ell\} P\{X(n(t-s)) \leq [nbt] - [nbs] + \ell - m\}.$$

Rearranging simplifies the sum to

$$\sum_{k<0} P\{X(n(t-s)) \leq [nbt] - [nbs] + k\}$$

$$\leq \sum_{k<0} P\{X(n(t-s)) \leq [nb(t-s)] + k + 1\}$$

(4.29)
$$= E[(X(n(t-s)) - [nb(t-s)] - 1)_-]$$

$$\leq C_0 \sqrt{n(t-s)} + C_1.$$

Above we first used $[nbt] - [nbs] \leq [nb(t-s)] + 1$. Assumption (2.6) gives uniform integrability to all moments of $u^{-1/2}(X(u) - [bu])$ as $u \to \infty$. Due to the $-1$ inside the expectation, we need the constant $C_1$ in (4.29), $C_0\sqrt{n(t-s)}$ alone will not do as $t - s \to 0$.

We leave the corresponding calculation for $P(B_{m,1})$ to the reader, and consider the lemma proved. $\square$

We are ready to prove Proposition 4.2. Apply (4.28) to each sum in (4.26), remove the squares with $2x^2 \leq x + x^3$, and let $C$ change its value from line

to line:

$$E[(\overline{Y}_n(t) - \overline{Y}_n(s))^6] \leq C\{\sqrt{n(t-s)} + 1 + (\sqrt{n(t-s)} + 1)^3\}$$
$$\leq C\{\sqrt{n(t-s)} + n^{3/2}(t-s)^{3/2} + 1\}.$$

Proposition 4.2 is proved. □

4.3. *Switch to discrete-time process.* The processes whose weak convergence is claimed in Theorem 2.2 are $n^{-1/4}Y_n(t)$. As observed earlier, it is equivalent to prove convergence for the centered processes $n^{-1/4}\overline{Y}_n(t)$. The moment estimate in (4.25) is not good enough for tightness, but we can get around this by a suitable time discretization. The forthcoming Lemma 4.11 contains an estimate that gives tightness for the time-discretized process we next define.

Fix two constants $\alpha, \beta > 0$ such that

(4.30) $$\tfrac{5}{4} + \alpha < \beta < \tfrac{3}{2}.$$

Let

(4.31) $$\widetilde{W}_n(t) = n^{-1/4}\overline{Y}_n(n^{-\beta}[n^\beta t])$$
$$= n^{-1/4}(Y_n(n^{-\beta}[n^\beta t]) - EY_n(n^{-\beta}[n^\beta t])).$$

In this section we show that it suffices to prove the weak convergence of process $\widetilde{W}_n$ by showing that $n^{-1/4}Y_n$ and $\widetilde{W}_n$ come uniformly close on compact time intervals.

PROPOSITION 4.3. *For any $0 < T < \infty$ and $\varepsilon > 0$,*

$$\lim_{n \to \infty} P\left\{\sup_{0 \leq t \leq T} |n^{-1/4}Y_n(t) - \widetilde{W}_n(t)| \geq \varepsilon\right\} = 0.$$

PROOF. Because the expectations $n^{-1/4}EY_n(n^{-\beta}[n^\beta t])$ vanish uniformly over $0 \leq t \leq T$ by Lemma 4.6, it suffices to prove

(4.32) $$\lim_{n \to \infty} P\left\{\sup_{0 \leq t \leq T} |Y_n(t) - Y_n(n^{-\beta}[n^\beta t])| \geq n^{1/4}\varepsilon\right\} = 0.$$

To prove (4.32), we consider the ways in which $Y_n(t)$ can differ from $Y_n(n^{-\beta}[n^\beta t])$ some time during $[0, T]$. First Lemma 4.9 shows that particles that start off at least at distance $n^{1/2+\alpha}$ from $[n\bar{y}]$ do not contribute to $Y_n(\cdot)$ during time interval $[0, T]$, in the $n \to \infty$ limit.

LEMMA 4.9. *Let*

$$N_1(T) = \sum_{m \leq [n\bar{y}] - n^{1/2+\alpha}} \sum_{j=1}^{\eta_0^n(m)} \mathbf{1}\{X_{m,j}^n(nt) \geq [n\bar{y}] + [nbt] \text{ for some } 0 \leq t \leq T\}$$
(4.33)
$$+ \sum_{m \geq [n\bar{y}] + n^{1/2+\alpha}} \sum_{j=1}^{\eta_0^n(m)} \mathbf{1}\{X_{m,j}^n(nt) \leq [n\bar{y}] + [nbt] \text{ for some } 0 \leq t \leq T\}$$

*be the number of particles that start at least at distance $n^{1/2+\alpha}$ from $[n\bar{y}]$, and reach the characteristic some time during $[0, nT]$. Then for a fixed $T$, $EN_1(T) \to 0$ as $n \to \infty$.*

PROOF. We handle the first sum in the definition of $N_1(T)$ and omit the similar argument for the other sum. Fix a positive integer $M$ large enough so that $1/2 - \alpha(2M-1) < 0$. Again $X(\cdot)$ denotes a random walk starting at the origin, and $C$ denotes a constant whose value may change from line to line. As $E\eta_0^n(m)$ is uniformly bounded [assumption (2.8)], and by an application of Doob's inequality to the martingale $X(t) - bt$, the expectation of the first sum in (4.33) is bounded by

$$C \sum_{\ell \geq n^{1/2+\alpha}} P\left\{\sup_{0 \leq t \leq T}(X(nt) - nbt) \geq \ell - 1\right\}$$
$$\leq C \sum_{\ell \geq n^{1/2+\alpha}} \ell^{-2M} E[(X(nT) - nbT)_+^{2M}]$$
$$\leq C \sum_{\ell \geq n^{1/2+\alpha}} \ell^{-2M} n^M \leq C n^{1/2 - \alpha(2M-1)}.$$

The expectation $E[(X(nT) - nbT)_+^{2M}]$ is $O(n^M)$, as suggested by the central limit theorem, due to uniform integrability guaranteed by assumption (2.6). □

Let

$$N_2 = \sum_{m\,:\,|m - [n\bar{y}]| \leq n^{1/2+\alpha}} \eta_0^n(m)$$

be the number of particles initially within distance $n^{1/2+\alpha}$ of $[n\bar{y}]$. Fix a constant $c$ so that

$$\lim_{n \to \infty} P\{N_2 \geq c n^{1/2+\alpha}\} = 0.$$

The event in (4.32) is contained in the event

$$\bigcup_{0 \leq k \leq [Tn^\beta]} \left\{\sup_{kn^{-\beta} \leq t \leq (k+1)n^{-\beta}} |Y_n(t) - Y_n(n^{-\beta} k)| \geq n^{1/4} \varepsilon\right\}.$$

For a fixed $k$, the event in braces implies that at least one of these two scenarios takes place:

(i) At least $\frac{1}{2}\varepsilon n^{1/4}$ particles cross the discretized characteristic $s \mapsto [n\bar{y}] + [bs]$ during time interval $s \in [n^{1-\beta}k, n^{1-\beta}(k+1)]$ by jumping. (Note that the time interval has been put in the microscopic time scale of the particles.) On the event $\{N_1(T) = 0\}$, these particles must be among the $N_2$ particles initially within $n^{1/2+\alpha}$ of $[n\bar{y}]$. Consequently, conditioned on $\{N_1(T) = 0\}$, the probability of this event is bounded by the probability that $N_2$ independent rate 1 random walks altogether experience at least $\frac{1}{2}\varepsilon n^{1/4}$ jumps in a time interval of length $n^{1-\beta}$.

(ii) At least $\frac{1}{2}\varepsilon n^{1/4}$ particles cross the discretized characteristic during time interval $[n^{1-\beta}k, n^{1-\beta}(k+1)]$ by staying put while the characteristic crosses the location of these particles. For large enough $n$, the distance between the endpoints $[n\bar{y}] + [n^{1-\beta}b(k+1)]$ and $[n\bar{y}] + [n^{1-\beta}bk]$ of the characteristic is at most 1. Hence, at most 1 site moves from one side of the characteristic to the other during this time interval, and so these $\frac{1}{2}\varepsilon n^{1/4}$ particles must sit on a unique site $x_k$ at time $n^{1-\beta}k$.

Accounting for all the possibilities gives the bound below. $\Pi(cn^{3/2+\alpha-\beta})$ is a mean $cn^{3/2+\alpha-\beta}$ Poisson random variable and represents the total number of jumps among $cn^{1/2+\alpha}$ independent particles during a time interval of length $n^{1-\beta}$,

$$P\left\{\sup_{0 \le t \le T} |Y_n(t) - Y_n(n^{-\beta}[n^\beta t])| \ge n^{1/4}\varepsilon\right\}$$

(4.34) $\quad \le P\{N_1(T) \ge 1\} + P\{N_2 \ge cn^{1/2+\alpha}\}$

$$+ \sum_{k=0}^{[Tn^\beta]} (P\{\Pi(cn^{3/2+\alpha-\beta}) \ge \tfrac{1}{2}n^{1/4}\varepsilon\} + P\{\eta_{kn^{1-\beta}}^n(x_k) \ge \tfrac{1}{2}n^{1/4}\varepsilon\}).$$

The probabilities $P\{N_1(T) \ge 1\}$ and $P\{N_2 \ge cn^{1/2+\alpha}\}$ vanish as $n \to \infty$ by Lemma 4.9 and choice of $c$. $\Pi(cn^{3/2+\alpha-\beta})$ is stochastically larger than a sum of $M_n = [cn^{3/2+\alpha-\beta}]$ i.i.d. mean 1 Poisson variables, and so a standard large deviation estimate gives

$$P\{\Pi(cn^{3/2+\alpha-\beta}) \ge \tfrac{1}{2}n^{1/4}\varepsilon\} \le \exp\{-M_n I(\tfrac{1}{2}M_n^{-1}n^{1/4}\varepsilon)\},$$

where $I$ is the Cramér rate function for the Poisson(1) distribution. By the choice of $\alpha$ and $\beta$, $M_n \ge n^\alpha$, while $M_n^{-1}n^{1/4} \to \infty$. Consequently, there are constants $0 < C_0, C_1 < \infty$,

$$\sum_{k=0}^{[Tn^\beta]} P\{\Pi(cn^{3/2+\alpha-\beta}) \ge \tfrac{1}{2}n^{1/4}\varepsilon\} \le C_0 n^\beta \exp(-C_1 n^\alpha) \to 0.$$

To treat the last term in (4.34), we derive a moment estimate for the occupation variables uniformly over space and time.

LEMMA 4.10. *Let $k \in \mathbf{N}$ and suppose initially $\sup_{m \in \mathbf{Z}} E[\eta_0(m)^k] < \infty$. Then*

$$\sup_{m \in \mathbf{Z}, t \geq 0} E[\eta_t(m)^k] < \infty.$$

PROOF. Fix $x \in \mathbf{Z}$ and $t > 0$, and let

$$\zeta_m = \sum_{j=1}^{\eta_0(m)} \mathbf{1}\{X_{m,j}(t) = x\}$$

be the number of particles initially at $m$ who find themselves at $x$ at time $t$. Then

$$E[\eta_t(x)^k] = E\left[\left(\sum_{m \in \mathbf{Z}} \zeta_m\right)^k\right]$$

$$= \sum_{m_1, m_2, \ldots, m_k \in \mathbf{Z}} E[\zeta_{m_1} \zeta_{m_2} \cdots \zeta_{m_k}]$$

$$= \sum_{b=1}^{k} \sum_{(m_1, m_2, \ldots, m_b)} \sum_{(k_1, k_2, \ldots, k_b)} \binom{k}{k_1 k_2 \cdots k_b} E[\zeta_{m_1}^{k_1}] E[\zeta_{m_2}^{k_2}] \cdots E[\zeta_{m_b}^{k_b}].$$

On the last line above we arrange the sum over all $k$-tuples $(m_1, \ldots, m_k) \in \mathbf{Z}^k$ according to the number $b$ of distinct sites among $m_1, \ldots, m_k$. The second sum on the last line is over $b$-tuples $(m_1, \ldots, m_b)$ of distinct sites from $\mathbf{Z}$. The third sum is over $b$-tuples $(k_1, \ldots, k_b)$ of positive integers such that $k_1 + \cdots + k_b = k$, and $\binom{k}{k_1 k_2 \cdots k_b}$ counts the number of ways $k_1$ $m_1$'s, $k_2$ $m_2$'s, and so on can be arranged into a $k$-tuple. Since $m_1, \ldots, m_b$ are distinct, $\zeta_{m_1}, \ldots, \zeta_{m_b}$ are independent.

Calculating as in the proof of Lemma 4.7 gives the bound

$$E[\zeta_m^k] \leq E[\eta_0(m)^k] p_t(m, x),$$

where $p_t(m, x) = p_t(0, x - m)$ is the translation-invariant transition probability of the underlying random walk of the particles. Substituting this back above gives the upper bound

$$\sum_{b=1}^{k} \sum_{(m_1, m_2, \ldots, m_b)} \sum_{(k_1, k_2, \ldots, k_b)} \binom{k}{k_1 k_2 \cdots k_b}$$

(4.35)

$$\times \left\{\prod_{i=1}^{b} p_t(0, x - m_i)\right\} \cdot \left\{\prod_{i=1}^{b} E[\eta_0(m_i)^{k_i}]\right\}.$$

Hölder's inequality, the moment assumption and $k_1 + \cdots + k_b = k$ give

$$\prod_{i=1}^{b} E[\eta_0(m_i)^{k_i}] \leq \prod_{i=1}^{b} E[\eta_0(m_i)^k]^{k_i/k} \leq C.$$

After this, sum the probabilities $p_t(0, x - m_i)$ over each index $m_i$ in (4.35). This leaves

$$\sum_{b=1}^{k} \sum_{(k_1, k_2, \ldots, k_b)} \binom{k}{k_1 k_2 \cdots k_b},$$

which is a constant that depends on $k$. □

We turn to the last term of (4.34),

$$\sum_{k=0}^{[Tn^\beta]} P\{\eta^n_{kn^{1-\beta}}(x_k) \geq \tfrac{1}{2} n^{1/4} \varepsilon\} \leq (Tn^\beta + 1) 2^6 \varepsilon^{-6} n^{-3/2} \sup_{x,t,n} E[\eta^n_t(x)^6].$$

Since $\sup_{x,t,n} E[\eta^n_t(x)^6] < \infty$ by the moment hypothesis (2.8) and Lemma 4.10, and $\beta - 3/2 < 0$, the right-hand side vanishes as $n \to \infty$.

We have shown that the right-hand side of the inequality in (4.34) vanishes as $n \to \infty$, and, thereby, proved Proposition 4.3. □

4.4. *Weak convergence.* We first verify tightness of the discrete-time processes. Let

$$\lambda = (3 - 2\beta)/(2\beta) \in (0, \tfrac{1}{5}).$$

LEMMA 4.11. *Fix $0 < T < \infty$. Then there exists a constant $C$ such that for all $0 \leq t_1 \leq t \leq t_2 \leq T$ and all $n$,*

(4.36) $\quad E[|\widetilde{W}_n(t) - \widetilde{W}_n(t_1)|^3 |\widetilde{W}_n(t_2) - \widetilde{W}_n(t)|^3] \leq C(t_2 - t_1)^{1+\lambda}.$

PROOF. If $t_2 - t_1 < n^{-\beta}$, then necessarily either $[n^\beta t] = [n^\beta t_1]$ or $[n^\beta t_2] = [n^\beta t]$. In either case, the left-hand side of (4.36) vanishes and the inequality holds trivially. So we may suppose $t_2 - t_1 \geq n^{-\beta}$.

By the Schwarz inequality and $2xy \leq x^2 + y^2$,

$$E[|\widetilde{W}_n(t) - \widetilde{W}_n(t_1)|^3 |\widetilde{W}_n(t_2) - \widetilde{W}_n(t)|^3]$$
$$\leq E[|\widetilde{W}_n(t) - \widetilde{W}_n(t_1)|^6]^{1/2} E[|\widetilde{W}_n(t_2) - \widetilde{W}_n(t)|^6]^{1/2}$$
$$\leq E[|\widetilde{W}_n(t) - \widetilde{W}_n(t_1)|^6] + E[|\widetilde{W}_n(t_2) - \widetilde{W}_n(t)|^6].$$

Apply (4.25) multiplied by $n^{-3/2}$ to both terms, ignoring the constant $C$ in the front, to get the upper bound

$$
\begin{aligned}
(4.37) \quad & \frac{1}{n}\left(\frac{[n^\beta t] - [n^\beta t_1]}{n^\beta}\right)^{1/2} + \frac{1}{n}\left(\frac{[n^\beta t_2] - [n^\beta t]}{n^\beta}\right)^{1/2} \\
& + \left(\frac{[n^\beta t] - [n^\beta t_1]}{n^\beta}\right)^{3/2} + \left(\frac{[n^\beta t_2] - [n^\beta t]}{n^\beta}\right)^{3/2} + 2n^{-3/2}.
\end{aligned}
$$

Since $t_2 - t_1 \geq n^{-\beta}$,

$$\frac{[n^\beta t] - [n^\beta t_1]}{n^\beta} \leq \frac{n^\beta t - n^\beta t_1 + 1}{n^\beta} = t - t_1 + \frac{1}{n^\beta} \leq 2(t_2 - t_1),$$

and also

$$\frac{1}{n} = (n^{-\beta})^{1/\beta} \leq (t_2 - t_1)^{1/\beta},$$

so the first term in (4.37) is bounded by

$$\frac{1}{n}\left(\frac{[n^\beta t] - [n^\beta t_1]}{n^\beta}\right)^{1/2} \leq 2^{1/2}(t_2 - t_1)^{1/2+1/\beta} \leq 2^{1/2} T^{(\beta-1)/(2\beta)}(t_2 - t_1)^{1+\lambda}.$$

Apply similar reasoning to the other terms. Note also that

$$(t_2 - t_1)^{3/2} \leq T^{1/2-\lambda}(t_2 - t_1)^{1+\lambda}$$

and

$$2n^{-3/2} = 2(n^{-\beta})^{1+\lambda} \leq 2(t_2 - t_1)^{1+\lambda}.$$

Collecting terms gives (4.36). □

Propositions 4.1 and 4.3 imply that the finite-dimensional distributions of the process $\widetilde{W}_n$ converge to those of $Z$ defined in Theorem 2.2. This and Lemma 4.11 are the hypotheses needed for Theorem 15.6 in [3]. We conclude that the processes $\widetilde{W}_n$ converge to the process $Z$ on the space $D_{\mathbf{R}}[0, T]$. Proposition 4.3 then implies that the processes $n^{-1/4} Y_n$ converge to $Z$ on the space $D_{\mathbf{R}}[0, T]$. Since $T$ is arbitrary, the convergence holds on $D_{\mathbf{R}}[0, \infty)$. This completes the proof of Theorem 2.2. Instead of Theorem 15.6 in [3], one can use Theorem 8.8 on page 139 of [6].

4.5. *Proof of Theorem* 2.3. Let $N_1(\bar{y}, T)$ denote the random variable defined by (4.33) to display its dependence on $\bar{y}$. On the event

$$(4.38) \qquad \{N_1(\bar{y}_1, T) = N_1(\bar{y}_2, T) = \cdots = N_1(\bar{y}_k, T) = 0\},$$

the processes $Y_n(\bar{y}_1, \cdot), Y_n(\bar{y}_2, \cdot), \ldots, Y_n(\bar{y}_k, \cdot)$, restricted to $[0, T]$, depend on disjoint collections of independent random walks if $n$ is large enough. The probability of (4.38) converges to 1 as $n \to \infty$ by Lemma 4.9. Consequently, the restrictions to $[0, T]$ of the processes $Y_n(\bar{y}_i, \cdot)$ become independent in the limit. To prove the tightness of the joint process $n^{-1/4}(Y_n(\bar{y}_1, \cdot), Y_n(\bar{y}_2, \cdot), \ldots, Y_n(\bar{y}_k, \cdot))$ on the space $D_{\mathbf{R}^k}[0, \infty)$, apply Theorem 8.8 from page 139 of [6] to the discrete-time process $(\widetilde{W}_n(\bar{y}_1, \cdot), \widetilde{W}_n(\bar{y}_2, \cdot), \ldots, \widetilde{W}_n(\bar{y}_k, \cdot))$, each component defined as in (4.31). The proof of Lemma 4.11 can be adapted to the multivariate case. We omit the details, and consider Theorem 2.3 proved.

4.6. *Proof of Theorem* 2.1. Limit (2.16) follows from (2.19) and Theorem 2.2. Subsequently, the hydrodynamic limit (2.15) follows from limit (2.16) and assumption (2.7).

To prove (2.17), we take a two-sided Brownian motion $B(\cdot)$, and create a bi-infinite sequence of random times

$$\cdots \leq T_{n,-2} \leq T_{n,-1} \leq 0 = T_{n,0} \leq T_{n,1} \leq T_{n,2} \leq \cdots$$

such that

(4.39) $$\sigma_0^n(x) = \sum_{m=1}^{x} \rho_0\left(\frac{m}{n}\right) + \sqrt{n} B(n^{-1} T_{n,x})$$

can be taken as the initial condition. (If $x < 0$, the sum actually ranges over $x + 1 \leq m \leq 0$.) To achieve this, apply the usual Skorokhod embedding (see, e.g., Section 7.6 in [5]) to the independent mean-zero random variables $\{\eta_0^n(x) - \rho_0(\frac{x}{n}) : x \in \mathbf{Z}\}$ and the two-sided Brownian motion $B_n(s) \equiv n^{1/2} B(\frac{s}{n})$. Embed $\{\eta_0^n(x) - \rho_0^n(x) : x > 0\}$ in the positive half of $B_n$, the remaining random variables in the negative half of $B_n$. Then the increments $\{T_{n,x} - T_{n,x-1} : x \in \mathbf{Z}\}$ are independent with means

(4.40) $$E(T_{n,x} - T_{n,x-1}) = \text{Var}[\eta_0^n(x)] = v_0\left(\frac{x}{n}\right),$$

and we have the equality in distribution of processes

$$\{B_n(T_{n,x}) - B_n(T_{n,x-1}) : x \in \mathbf{Z}\} \stackrel{d}{=} \left\{\eta_0^n(x) - \rho_0\left(\frac{x}{n}\right) : x \in \mathbf{Z}\right\}.$$

Note that we have been using (2.11) which is assumed for this part of Theorem 2.1. Now it is clear that $\sigma_0^n$, defined by (4.39), has the right distribution to serve as the initial height function.

Next we observe the central limit theorem for $\sigma_0^n$. From (4.39), for $y \in \mathbf{R}$,

(4.41) $$\frac{\sigma_0^n([ny]) - n u_0(y)}{\sqrt{n}} = B(n^{-1} T_{n,[ny]})$$
$$+ \sqrt{n} \int_0^y \left(\rho_0\left(\frac{[ns]}{n}\right) - \rho_0(s)\right) ds + O(n^{-1/2}).$$

One can show that

$$\lim_{n \to \infty} n^{-1} T_{n,[ny]} = \int_0^y v_0(s) \, ds \qquad \text{in probability.}$$

This follows from Chebyshev's inequality and the moment bound

$$E[(T_{n,m} - T_{n,m-1})^2] \le CE[(\eta_0^n(m) - E\eta_0^n(m))^4] \le C_1.$$

As $n \to \infty$, the integral term in (4.41) vanishes by the Hölder property (2.10) of $\rho_0$. Consequently,

$$\frac{\sigma_0^n([ny]) - nu_0(y)}{\sqrt{n}} \to B\left(\int_0^y v_0(s) \, ds\right)$$

in probability. Finally, (2.17) follows from this and (2.16).

**5. Proof for Hammersley's process.** The proof of Theorem 3.2 is based on the approach and estimates derived in [18]. To save space, we refer to that paper for all the groundwork.

We construct the initial configuration by Skorohod's representation, so that

$$z_0^n(i) = nu_0(i/n) + n^{1/2} B(n^{-1} T_{n,i}),$$

where for each fixed $n$, $\{T_{n,i} : i \in \mathbf{Z}\}$ are the hitting times of appropriate random intervals independent of the two-sided Brownian motion $B(\cdot)$. Section 8 in [18] discusses this construction. For $(x,t) \in \mathbf{R} \times [0,\infty)$, set

$$\zeta_t^n(x) = n^{-1/2} \{z_{nt}^n([nx]) - nu(x,t)\}.$$

In particular, in terms of the Brownian motion at $t = 0$,

(5.1) $\qquad \zeta_0^n(y) = B(n^{-1} T_{n,[ny]}) + n^{1/2}(u_0([ny]/n) - u_0(y)).$

Hammersley's process has a special graphical construction in terms of increasing sequences among rate 1 Poisson points on the space-time plane. This representation goes back to [9] and [1]. It can be expressed as follows:

(5.2) $\qquad\qquad z_t^n(k) = \inf_{i \, : \, i \le k} \{z_0^n(i) + \Gamma_t^{n,i}(k-i)\},$

where $\Gamma_t^{n,i}(m)$ is the minimal positive $h$ such that the space-time rectangle $(z_0^n(i), z_0^n(i) + h] \times (0, t]$ contains at least $m$ Poisson points on an increasing path. This appears as equation (45) in [18]. But observe that compared to [18], in the present paper time arguments have become subscripts and space indices have become arguments in parentheses.

5.1. *Upper bound.* Following Section 6.1.4 on page 172 in [18] gives $Y_n \leq R_{n,3} + C_0$ where $C_0$ is a constant and

$$R_{n,3} = \sup_{y \in I(x,t)} \left\{ \Gamma_{nt}^{n,[ny]}([nx] - [ny]) - \frac{([nx] - [ny])^2}{4nt} \right\}.$$

By Lemmas 4.4 and 4.5 in [18], $P(R_{n,3} \geq Cn^{1/3} \log n) \to 0$ if $C$ is fixed large enough.

5.2. *Lower bound.* Let $\varepsilon > 0$. We shall show that for some constant $0 < C < \infty$, $Y_n \geq -Cn^{1/3} \log n$ with probability at least $1 - \varepsilon$ for large $n$.

Let $i_n$ be the minimal microscopic minimizer for $z_{nt}^n([nx])$ defined by (59) in [18]. By Lemma 5.3 in [18], $\text{dist}(n^{-1}i_n, I(x,t)) \to 0$ in probability. Let $\delta > 0$ be as in Assumption E. Then for large enough $n$, there exists a random $y_n \in I(x,t)$ such that

(5.3) $$P\{|n^{-1}i_n - y_n| < \delta\} \geq 1 - \varepsilon/4.$$

Assume now we are on the event in braces in (5.3). Let $[a_1, a_2]$ be a compact interval that contains the $\delta$-neighborhood of $I(x,t)$. Following the calculation in Section 6.1.2 on page 170 in [18] gives

$$Y_n \geq \left\{ \Gamma_{nt}^{n,i_n}([nx] - i_n) - \frac{(nx - i_n)^2}{4nt} \right\} + n^{1/2} \{\zeta_0^n(n^{-1}i_n) - \zeta_0^n(y_n)\}$$
$$+ n\{\Phi(n^{-1}i_n) - \Phi(y_n)\}$$
$$\geq R_{n,1} + n^{1/2} \{\zeta_0^n(n^{-1}i_n) - \zeta_0^n(y_n)\} + c_1 n^{-1}(i_n - ny_n)^2.$$

We used the definition

$$R_{n,1} = \min\left\{ \Gamma_{nt}^{n,i}([nx] - i) - \frac{(nx - i)^2}{4nt} : na_1 \leq i \leq na_2 \right\},$$

and applied assumption (3.6). By Lemma 4.3 in [18], for large enough $n$,

(5.4) $$P\{R_{n,1} \geq -Cn^{1/3} \log n\} \geq 1 - \varepsilon/8$$

if $C$ is fixed large enough. Of the lower bound of $Y_n$ given above, it remains to handle this part:

(5.5)
$$n^{1/2} \{\zeta_0^n(n^{-1}i_n) - \zeta_0^n(y_n)\} + c_1 n^{-1}(i_n - ny_n)^2$$
$$\geq n^{1/2} \{B(n^{-1}T_{n,i_n}) - B(n^{-1}T_{n,[ny_n]})\} + c_1 n^{-1}(i_n - ny_n)^2 - c_2.$$

The last constant $-c_2$ accounts for the last terms in (5.1), actually of order $n^{-1/2}$.

Let $[b_1, b_2]$ be an interval such that $\int_0^y v_0(s)\,ds \in [b_1 + 1, b_2 - 1]$ for all $y \in [a_1, a_2]$. Let

$$w(u) = \sup\{|B(s) - B(r)| : s, r \in [b_1, b_2], |s - r| \leq u\}$$

be the modulus of continuity of Brownian motion on the interval $[b_1, b_2]$ of length $A = b_2 - b_1$. By Lévy's theorem,

$$\limsup_{u \to 0}(2u \log(A/u))^{-1/2} w(u) = 1 \quad \text{a.s.}$$

Choose $\delta_0 \in (0, 1/2)$ so that

(5.6) $\quad P\{w(u) \leq 2(u \log(A/u))^{1/2} \text{ for all } u \leq \delta_0\} \geq 1 - \varepsilon/16.$

The following law of large numbers holds: for all large enough $n$,

(5.7) $\quad P\left\{\sup_{y \in [a_1, a_2]} \left| n^{-1} T_{n,[ny]} - \int_0^y v_0(s)\, ds \right| \leq \delta_0 \right\} \geq 1 - \varepsilon/32.$

This can be proved as Lemma 8.1 in [18] is proved. The fourth moment assumption enters here, because the hitting times satisfy

$$E[(T_{n,i} - T_{n,i-1})^2] \leq CE[\eta_0^n(i)^4].$$

On the event (5.3), both $y_n$ and $n^{-1} i_n$ lie in $[a_1, a_2]$. Then on the event in (5.7), both $n^{-1} T_{n,i_n}$ and $n^{-1} T_{n,[ny_n]}$ lie in $[b_1, b_2]$ and so fall within the range of the event in (5.6).

Thus, on the intersection of the events in (5.3), (5.6) and (5.7), the right-hand side of the inequality in (5.5) is bounded below by

(5.8) $\quad -2|T_{n,i_n} - T_{n,[ny_n]}|^{1/2}(\log An)^{1/2} + c_1 n^{-1}(i_n - ny_n)^2 - c_3.$

The new constant $-c_3$ includes the earlier constant $-c_2$ from (5.5) and the case where

$$|T_{n,i_n} - T_{n,[ny_n]}| < 1,$$

as otherwise this factor should be in the denominator inside the logarithm. But since $u \log(A/u) \to 0$ as $u \to 0$, the case of small $|T_{n,i_n} - T_{n,[ny_n]}|$ can be bounded by a single constant.

The probabilities of the complements of the events in (5.3), (5.4), (5.6) and (5.7) add up to less than $\varepsilon/2$. On the intersection of these events, $Y_n$ is bounded below by $-Cn^{1/3} \log n$, plus the expression in (5.8). It remains to show that, with probability at least $1 - \varepsilon/2$, the expression in (5.8) is bounded below by $-Cn^{1/3} \log n$ for yet another constant $C$.

We begin with a simple general fact. Suppose $\{X_i\}$ are nonnegative, independent random variables with bounded means and variances. Let $0 < a, b < \infty$, and $C \geq 1 + EX_i$ and $C_1 \geq \text{Var}\, X_i$ for all $i$. Then,

(5.9) $\quad P\left\{\sup_{1 \leq \ell \leq an, n^{2/3} \leq k \leq bn} \sum_{i=\ell}^{\ell+k-1} X_i \geq 3Ck \right\} \leq C_1(a+b) n^{-1/3}.$

To prove this, note first that the event in braces implies that

$$\sum_{mn^{2/3} \leq i < (m+1)n^{2/3}} X_i \geq Cn^{2/3} \qquad \text{for some } 1 \leq m \leq (a+b)n^{1/3}.$$

Then by Chebyshev's inequality,

$$P\bigg\{\sup_{1 \leq m \leq (a+b)n^{1/3}} \sum_{mn^{2/3} \leq i < (m+1)n^{2/3}} X_i \geq Cn^{2/3}\bigg\}$$

$$\leq P\bigg\{\sup_{1 \leq m \leq (a+b)n^{1/3}} \sum_{mn^{2/3} \leq i < (m+1)n^{2/3}} (X_i - EX_i) \geq n^{2/3}\bigg\}$$

$$\leq (a+b)n^{1/3}\frac{C_1 n^{2/3}}{n^{4/3}}.$$

Applying (5.9) to $|T_{n,i_n} - T_{n,[ny_n]}|$ (these differences are sums of independent, nonnegative increments), we conclude that on the event (5.3), for large enough $n$ with probability at least $1 - \varepsilon/2$,

$$|T_{n,i_n} - T_{n,[ny_n]}| \leq C(|i_n - ny_n| \vee n^{2/3}).$$

From this, the expression in (5.8) is bounded below by

$$(5.10) \quad -2C(|i_n - ny_n|^{1/2} \vee n^{1/3})(\log An)^{1/2} + c_1 n^{-1}(i_n - ny_n)^2 - c_3.$$

This last expression is no less than $-Cn^{1/3}(\log n)^{2/3}$ for a suitable (new) constant $0 < C < \infty$. This a lower bound for $Y_n$ with probability at least $1 - \varepsilon$.

DEPARTMENT OF MATHEMATICS
UNIVERSITY OF WISCONSIN–MADISON
MADISON, WISCONSIN 53706-1388
USA
E-MAIL: seppalai@math.wisc.edu